\documentclass[11pt,twoside]{article}
\usepackage{epsfig,rotating,latexsym,makeidx,amscd,wrapfig,amssymb,color}
\usepackage{amsmath,mathrsfs,amsfonts,stmaryrd}
\title{A numerical method for computing the Jordan Canonical Form
}
\author{
Zhonggang Zeng\thanks{Department of Mathematics, Northeastern
Illinois University, Chicago, IL 60625 (email: {\tt zzeng@neiu.edu}).
Research supported in part by NSF under Grant DMS-0412003.}
\and
Tien-Yien Li\thanks{Department of Mathematics, Michigan State
University, East Lansing, MI 48824 (email: {\tt li@math.msu.edu}).
Research supported in part by NSF under Grant DMS-0411165.}
}
\DeclareSymbolFont{AMSb}{U}{msb}{m}{n}
\DeclareMathSymbol{\bdC}{\mathbin}{AMSb}{'103}
\DeclareMathSymbol{\bdR}{\mathbin}{AMSb}{'122}
\DeclareMathSymbol{\N}{\mathbin}{AMSb}{'116}
\DeclareMathAlphabet{\mathpzc}{OT1}{pzc}{m}{it}
\newcommand{\h}{{{\mbox{\tiny $\mathsf{H}$}}}}
\newcommand{\cR}[1]{{\cal R}\left(\,#1\,\right)}
\newcommand{\cN}[1]{{\cal K}\left(\,#1\,\right)}

\newcommand{\cF}{{\cal F}}

\newcommand{\cU}{{\cal U}}
\newcommand{\cV}{{\cal V}}
\newcommand{\cW}{{\cal W}}

\newcommand{\blb}{\big[\,}
\newcommand{\brb}{\, \big]}

\newcommand{\spn}[1]{\mathpzc{span}\big\{\,#1\,\big\}}
\newcommand{\eig}[1]{\Lambda \big(\,#1\,\big)}
\newcommand{\dist}[1]{\mathpzc{dist}\big(\,#1\,\big)}
\newcommand{\diag}[1]{\mathpzc{diag}\big(\,#1\,\big)}
\newcommand{\codim}[1]{\mathpzc{codim}\big(\,#1\,\big)}
\newcommand{\dm}[1]{\mathpzc{dim}\big(\,#1\,\big)}
\newcommand{\rank}[1]{\mathpzc{rank}\big(\,#1\,\big)}
\newcommand{\ranka}[2]{\mathpzc{rank}_{#1}\big(\,#2\,\big)}
\newcommand{\nullity}[1]{\mathpzc{nullity}\big(\,#1\,\big)}

\newcommand{\bdf}{\mathbf{f}}
\newcommand{\bdg}{\mathbf{g}}
\newcommand{\bdu}{\mathbf{u}}
\newcommand{\bdv}{\mathbf{v}}

\newcommand{\bdx}{\mathbf{x}}
\newcommand{\bdy}{\mathbf{y}}
\newcommand{\bdz}{\mathbf{z}}

\newcommand{\bdb}{\mathbf{b}}
\newcommand{\bdc}{\mathbf{c}}

\newcommand{\bde}{\mathbf{e}}
\newcommand{\bdh}{\mathbf{h}}

\newcommand{\bdp}{\mathbf{p}}
\newcommand{\bdq}{\mathbf{q}}

\newcommand{\bds}{\mathbf{s}}

\newcommand{\bdo}{\mathbf{0}}
\newcommand{\la}{\lambda}
\newcommand{\al}{\alpha}

\newcommand{\dl}{\delta}

\newcommand{\eps}{\varepsilon}
\newcommand{\sg}{\sigma}

\newcommand{\hla}{\hat{\lambda}}
\newcommand{\hS}{\hat{S}}

\newcommand{\hU}{\hat{U}}
\newcommand{\divs}{\,\big|\,}
\newcommand{\mns}{\mbox{\footnotesize -}}
\newcommand{\pls}{\mbox{\raisebox{-.4ex}{\tiny $^+$\hspace{-0.3mm}}}}
\newcommand{\tms}{\mbox{\raisebox{-.4ex}{\tiny $^\times$\hspace{-0.4mm}}}}
\newcommand{\lmns}{\mbox{-}}
\newcommand{\lpls}{\mbox{\footnotesize +}}
\newcommand{\ltms}{\mbox{\hspace{-0.4mm}\footnotesize $\times$}}

\newcommand{\bnorm}[1]{\Big\|\, #1 \,\Big\|_2}

\newcommand{\prf}{\noindent {\bf Proof. \ }}
\newtheorem{defn}{Definition}
\newtheorem{lemma}{Lemma}
\newtheorem{theorem}{Theorem}
\newtheorem{prop}[theorem]{Proposition}

\newtheorem{example}{Example}

\newcommand{\qed}{${~} $ \hfill \raisebox{-0.3ex}{\LARGE $\Box$}}

\begin{document}

\thispagestyle{empty}
\begin{center}
{{\bf \huge A Numerical Method}
\linebreak
\linebreak {\Large for}
\linebreak
\linebreak ${ }_{ }$ \hspace{5mm} {\bf \huge Computing the Jordan Canonical
Form}}
\vspace{4mm}
\end{center}
\centerline{{\bf \Large Zhonggang Zeng} \ and \ {\Large \bf Tien-Yien Li}}
\vspace{2mm}
%\centerline{\bf \large Northeastern Illinois University}

\pagenumbering{roman}
\parskip-0.3mm
\message{Table of contents}
\addcontentsline{toc}{section}{\protect\numberline{}{}}{}
\tableofcontents

\newpage
\thispagestyle{empty}

\include{symbols}

\newpage
\setcounter{page}{1}

\parskip4mm

\pagenumbering{arabic}

\maketitle

\begin{abstract}
The Jordan Canonical Form of a matrix is highly sensitive to
perturbations, and its numerical computation remains a formidable challenge.
~This paper presents a regularization theory that establishes
a well-posed least squares problem of finding the nearest staircase
decomposition in the matrix bundle of the highest codimension.
~A two-staged algorithm is developed for computing the numerical
Jordan Canonical Form.
~At the first stage, the method calculates the Jordan structure of
the matrix and an initial approximation to the multiple eigenvalues.
~The staircase decomposition is then constructed
by an iterative algorithm at the second stage.
As a result, the numerical Jordan Canonical decomposition along with multiple
eigenvalues can be computed with high accuracy even if the underlying
matrix is perturbed.
\end{abstract}

{\bf keywords} ~~Jordan canonical form, eigenvalue, staircase form,

%\begin{AMS} 65F15 \end{AMS}

%\thispagestyle{plain}

\vspace{-3mm}
\section{Introduction}
\vspace{-5mm}

This paper presents an algorithm and a regularization theory for 
computing the Jordan Canonical Form accurately even if the matrix 
is perturbed.

The existence of the Jordan Canonical Form is one of the fundamental
theorems in linear algebra as an indispensable tool in
matrix theory and beyond.
~In practical applications, however, it is well documented that the Jordan
Canonical Form is extremely difficult, if not impossible, for numerical
computation
\cite[p.25]{eigtemp},
\cite[p.52]{bar-cam},
\cite[p.189]{chatelin-fraysse},
\cite[p.165]{chatelin},
\cite[p.146]{dem-book},
\cite[p.371]{gvl}, \cite[p.132]{hj}, \cite[p.22]{stew2}.
~In short, as remarked in a celebrated survey article by
Moler and Van Loan \cite{mvl}:
~``The difficulty is that the JCF cannot be computed using floating
point arithmetic. A single rounding error may cause some multiple
eigenvalue to become distinct or vise versa, altering the entire
structure of ~$J$~ and ~$P$.''

Indeed, defective multiple eigenvalues in a
non-trivial Jordan Canonical Form degrade to clusters of simple eigenvalues
in practical numerical computation.
~A main theme of the early attempts for
numerical computation of the Jordan Canonical Form is to locate a multiple
eigenvalue as the mean of a cluster that is selected from eigenvalues
computed by QR algorithm and, when it succeeds, the Jordan structure may be 
determined by computing a staircase form at the multiple eigenvalue.
This approach includes works of Kublanovskaya \cite{kubl} (1966),
Ruhe \cite{ruhe-70-bit} (1970), Sdridhar et al \cite{sri-jor} (1973),
and culminated in Golub and Wilkinson's review
\cite{golub-wilkinson} (1976) as well as K{\aa}gstr\"{o}m and Ruhe's
{\sc JNF} \cite{kagstrom-ruhe-jnf,kagstrom-ruhe} (1980).
~Theoretical issues have been analyzed in, e.g.
\cite{dem-thesis,demmel-86,wilk-84,wilk-86}.

However, the absence of a reliable method for identifying the proper
cluster renders a major difficulty for this approach.
~Even if the correct cluster can be identified, its arithmetic mean may not
be sufficiently accurate for identifying the Jordan structure, as shown in
Example \ref{e:sens} (\S \ref{sec:sens}).
~While improvements have been made steadily 
\cite{beelen-vandooren,lippert-edelman},
a qualitative approach is proposed in \cite{chatelin-fraysse}, and a partial
canonical form computation is studied in \cite{kag-wib},
``attempts to compute the Jordan canonical form of a matrix have not
been very successful'' as commented by Stewart in \cite[p. 22]{stew2}.

A related development is to find a well-conditioned matrix ~$G$~ such that
~$G^{\mns 1} A G$~ is block diagonal \cite[\S 7.6.3]{gvl}.
~Gu proved this approach is NP-hard \cite{gu95}, with a suggestion
that ``it is still possible that there are algorithms that can
solve most practical cases'' for the problem.
Another closely related problem is the computation of the Kronecker Canonical
Form for a matrix pencil ~$A-\la B$~ 
(see \cite{demmel-kagstrom,eek1,eek2,bk00}).
~For a given Jordan structure, a minimization method is proposed in
\cite{lip-edel-00} to find the nearest matrix with the same Jordan structure.

Multiple eigenvalues are multiple roots of the characteristic polynomial
of the underlying matrix.
~There is a perceived barrier of ``attainable accuracy'' associated
with multiple zeros of algebraic equations
which, in terms of number of digits,
is the larger one between data error and machine precision
divided by the multiplicity \cite{victorpan97,ypma,zeng-05}.
~Thus, as mentioned above, accurate computation of multiple
eigenvalues remains a major obstacle of computing the Jordan
Canonical Form using floating point arithmetic.
~Recently, a substantial progress has been achieved in computing
multiple roots of polynomials.
~An algorithm is developed in \cite{zeng-05} along
with a software package \cite{zeng_multroot}
that consistently determines multiple roots and their
multiplicity structures of a polynomial with remarkable accuracy
without using multiprecision
arithmetic even if the polynomial is perturbed.
~The method and results realized Kahan's observation in 1972 that multiple
roots are well behaved under perturbation when the multiplicity structure
is preserved \cite{kahan72}.

Similar to the methodology in \cite{zeng-05},
we propose a two-stage algorithm in this paper for computing
the numerical Jordan Canonical Form.
~To begin, we first find the Jordan structure in terms of
the Segre/Weyr characteristics at each distinct eigenvalue.
~With this structure as a constraint, the problem of computing the Jordan
Canonical Form is reformulated as a least squares problem.
~We then iteratively determine the accurate eigenvalues and a
{\em staircase decomposition}, and the Jordan decomposition can follow
as an option.

We must emphasize the {\em numerical} aspect of our algorithm that focuses on
computing the {\em numerical} Jordan Canonical Form of inexact matrices.
~The exact Jordan Canonical Form of a matrix with exact data
may be obtainable in many cases using symbolic computation
(see, e.g. \cite{olaz,fort-gia,giesbrecht,lzw}).
~Due to ill-posedness of the Jordan Canonical Form, however,
symbolic computation may not be suitable for applications where matrices
will most likely be perturbed in practice.
~For those applications, we must formulate the notion of the numerical
Jordan Canonical Form that is structurally invariant under small data
perturbation, and continuous in a neighborhood of the matrix with the exact
Jordan Canonical Form in question.

More precisely, matrices sharing a particular Jordan structure form
a matrix {\em bundle}, or, a manifold.
~For a given matrix ~$A$, ~we compute the exact Jordan Canonical
Form of the nearest matrix ~$\tilde{A}$~ in a bundle ~$\Pi$
~of the highest co-dimension within a neighborhood of ~$A$.
~Under this formulation, computing the numerical Jordan Canonical Form
of ~$A$~ should be a well-posed problem when ~$A$~ is sufficiently close
to bundle ~$\Pi$.
~In other words, under perturbation of sufficiently small magnitudes,
the deviation
of the numerical Jordan Canonical Form is tiny with the structure intact.

The main results of this paper can be summarized as follows.
~In \S \ref{s:reg}, we formulate a system of quadratic equations that uniquely
determines a local staircase decomposition at a multiple eigenvalue from a given
Jordan structure.
~Regularity theorems (Theorem \ref{unisc} and Theorem \ref{t:cplxreg}) in this
section establish the well-posedness of the staircase decomposition that ensures
accurate computation of multiple eigenvalues.
~Based on this regularity, the numerical unitary-staircase eigentriplets
is formulated in \S \ref{sec:sens}, along with the backward error measurement
and a proposed condition number.

In \S \ref{s:ctrip}, we present an iterative algorithm for computing the
well-posed unitary-staircase eigentriplet assuming the Jordan structure is
given.
~The algorithm employs the Gauss-Newton iteration whose local convergence
is a result of the regularity theorems given in \S \ref{s:reg}.
~The method itself can be used as a stand-alone algorithm for calculating the
nearest staircase/Jordan decomposition of a given structure, as demonstrated
via numerical examples in \S \ref{sec:numres1}.

The algorithm in \S \ref{s:ctrip} requires a priori knowledge of the
Jordan structure, which can be computed by an algorithm we propose in
\S \ref{s:compstru}.
~The algorithm employs a special purpose Hessenberg reduction and a
rank-revealing mechanism that produces the sequence of minimal
polynomials.
~Critically important in our algorithm is the application of the recently
established robust multiple root algorithm \cite{zeng-05} to those numerically
computed minimal polynomials in determining the Jordan structure as well as
an initial approximation of the multiple eigenvalues, providing
the crucial input items needed in the staircase algorithm in \S \ref{s:ctrip}.
~In \S \ref{s:overall}, we summarize the overall algorithm and
present numerical results%
%\footnote{The Matlab code of our algorithm and test scripts for all the 
%numerical examples can be 
%accessed at {\tt http://www.neiu.edu/$\sim$zzeng/numjcf.htm}
%}.

\vspace{-5mm}
\section{Preliminaries}
\vspace{-5mm}

\subsection{Notation and terminology}
\vspace{-5mm}

Throughout this paper, \label{'matrix'} matrices are denoted by upper case 
letters ~$A$, ~$B$, ~etc., ~and ~$O$~ 
denotes a zero matrix \label{'0matrix'} with known dimensions.
~Vectors \label{'vector'} are in columns and represented by lower case 
boldface letters like ~$\bdu$, ~$\bdv$~ and ~$\bdx$.
~A zero vector \label{'0vector'} is denoted by ~$\bdo$, ~or ~$\bdo_n$~ to 
emphasize the dimension.
~The notation ~$(\cdot)^\top$~ \label{'top'}
represents the transpose of a matrix or a vector ~$(\cdot)$,
~and ~$(\cdot)^\h$ ~is its Hermitian adjoint (or conjugate transpose).
~The fields of real and complex numbers are denoted by ~$\bdR$~ and
~$\bdC$~ \label{'RC'}respectively.

For any matrix ~$B \in \bdC^{m\tms n}$, \label{'rank'}
~the rank, nullity, range and kernel of ~$B$~ are denoted by
~$\rank{B}$, ~$\nullity{B}$, ~$\cR{B}$~ and ~$\cN{B}$~ respectively.
~The ~$n \times n$~ identity matrix \label{'idmatrix'} is ~$I_{n}$,
~or simply ~$I$~ when its size is clear.
~The column vectors of ~$I$~ are {\bf canonical vectors}
\label{'canonvec'} ~$\bde_1,\cdots,\bde_n$.
~A matrix ~$U \in \bdC^{m\tms n} $~ is said to be
{\bf unitary}\index{matrix!unitary} if ~$U^\h U = I$.
~A matrix ~$V \in \bdC^{m\tms (n\mns m)} $~ is called a unitary complement
of unitary matrix ~$U$~ if ~$[U,V]$~ is a square unitary matrix.
~Subspaces \label{'subspace'} of ~$\bdC^n$~ are denoted by calligraphic letters 
~${\cal X}$, ~${\cal Y}$, ~with dimensions \label{'dim'} ~$\dm{\cal X}$, 
~$\dm{\cal Y}$, ~etc., and ~${\cal X}^\perp$~ stands for the orthogonal
complement \label{'perp'} of subspace ~${\cal X}$.
~The set of distinct eigenvalues of ~$A$~ is the
{\bf spectrum}\label{'spectrum'}
of ~$A$~ and is denoted by
~$\eig{A}$\index{$\eig{\cdot}$: \ \ spectrum of $(\cdot)$}.

\vspace{-5mm}
\subsection{Segre and Weyr characteristics}
\vspace{-5mm}

The Jordan structure and the corresponding staircase structure 
(see \S \ref{s:scf}) 
of an eigenvalue can be characterized by the Segre characteristic
and the Weyr characteristic respectively.
~These two characteristics are conjugate partitions
of the algebraic multiplicity of the underlying eigenvalue.
~Here, a sequence of nonnegative integers ~$\{ k_1 \geq k_2 \geq \cdots \}$~
is called a
{\bf partition}\index{partition!of a positive integer}
of a positive integer ~$k$~ if ~$k = k_1+k_2+\cdots$.
~For such a partition, sequence
~$l_j = \max \Big\{\,i\,\Big|\, k_i \geq j \Big\}$, ~$j=1, 2, \cdots$
~is called the {\bf conjugate partition} of ~$\{k_1,k_2,\cdots \}$.
~For example, ~$[3,2,2,1]$~ is a partition of $8$ with
conjugate ~$[4,3,1]$~ and vice versa.

Let ~$\la$~ be an eigenvalue of ~$A$~ with an algebraic multiplicity ~$m$~
corresponding to elementary Jordan blocks of
orders ~$n_1\geq n_2 \geq \cdots \geq n_l > 0$.
~The {\em infinite} sequence
~$\{ n_1, \cdots, n_l, 0, 0, \cdots\}$~ is called
the {\bf Segre characteristic} of ~$A$~ associated with ~$\la$.
~The Segre characteristic forms a partition of the algebraic
multiplicity ~$m$~ of ~$\la$. ~Its conjugate partition is called
the {\bf Weyr characteristic} of ~$A$~ associated with ~$\la$.
~We also take the Weyr characteristic as an infinite sequence
for convenience.
~The nonzero part of such sequences will be called the {\em nonzero}
Segre/Weyr characteristics.
~Let ~$A$~ be an ~$n \times n$~ matrix with Weyr characteristic
~$\{m_1 \geq m_2 \geq \cdots\}$~ associated with an eigenvalue ~$\la$.
~Then \cite[Definition 3.6 and Lemma 3.2]{dem-edel}, for ~$j = 1, 2, \cdots$, 
\begin{eqnarray*}
%m_1 & = & \nullity{A-\la I} \;\; \equiv \;\; n-\rank{A-\la I}, \\
m_j & = &
\nullity{(A-\la I)^j} - \nullity{(A-\la I)^{j\mns 1}} 
\end{eqnarray*}
which immediately implies the uniqueness of
the two characteristics and their invariance under unitary
similarity transformations,
since the rank of ~$(PAP^{\mns 1}-\la I)^j = P(A-\la I)^j P^{\mns 1}$~ is
the same as the rank of ~$(A - \la I)^j$~ for ~$j = 1, 2, \cdots$.
~In particular, both characteristics are invariant under
Hessenberg reduction \cite[\S 7.4.3]{gvl}.

\vspace{-5mm}
\subsection{The staircase form} \label{s:scf}
\vspace{-5mm}
Discovered by Kublanovskaya \cite{kubl},
a matrix is associated with a
{\bf staircase form} given below.

\vspace{-4mm}
\begin{lemma} \label{lem:scf}
~~Let ~$A\in \bdC^{n\tms n}$~ be a matrix with nonzero Weyr characteristic
~$\left\{\,m_j\,\right\}_{j=1}^k$~ associated with an ~$m$-fold eigenvalue
~$\la$.
~For consecutive ~$j = 1, \cdots, k$, ~let ~$Y_j \in \bdC^{n\tms m_j}$~ be a
matrix satisfying
~$\cR{\blb Y_1, \cdots, Y_j \brb} =
\cN{ (A-\la I)^j }$.
~Then ~$\blb Y_1,\cdots,Y_k \brb$~ is of full rank and
\begin{eqnarray} 
\label{auus0}
A\blb Y_1,\cdots,Y_k \brb & = & \blb Y_1,\cdots,Y_k \brb
\left(\la I_{m} + S \right)  \\
& & \nonumber \\  \label{sl}
 \mbox{ \ where \ }
S & = & \mbox{\scriptsize $
\begin{array}{cl}
 \begin{array}{clll} \mbox{\small $m_1$} & \mbox{\small $m_{2}$} &
\cdots \, & \mbox{\small $m_k$} \;\;\;\;\;\; \end{array} & \\
\left[ \begin{array}{cccl}
\,O\, & S_{12} &  \cdots & S_{1k} \\
      & \ddots & \ddots & \vdots \\
      &        & \ddots & S_{k\mns 1,k} \\
      &        &        & \,O\, \end{array} \right]
&\begin{array}{l} \mbox{\hspace{-4mm}\raisebox{1pt}{\small $m_1$}} \\
\mbox{\hspace{-3mm}$\vdots$} \\
\mbox{\hspace{-4mm}\raisebox{-5pt}{\small $m_{k\mns 1}$}}
 \\ \mbox{\hspace{-4mm}\raisebox{-4pt}{\small $m_k$}} \end{array}
\end{array} $}
\end{eqnarray}
Furthermore, ~all super-diagonal blocks ~$S_{12},S_{23},\cdots,S_{k\mns 1,k}$~
are matrices of full rank.
\end{lemma}
\vspace{-4mm}

\prf
~Equation (\ref{auus0}) and the existence of 
~$S$ ~in (\ref{sl}) can be proved by a straightforward verfication 
using ~$\cR{Y_l} = \cN{(A-\la I)^l}$ ~for ~$l=1,\cdots, k$.
~From (\ref{sl}), we have
\[ (A-\la I)^{l\mns 1}\, Y_l \;=\; (A-\la I)^{l\mns 2}\,
\Big( Y_1 S_{1,l} + \cdots + Y_{l\mns 1} S_{l\mns 1,l} \Big)
\;=\; (A-\la I)^{l\mns 2}  Y_{l\mns 1} S_{l\mns 1,l}.
\]
This implies ~$S_{l\mns 1,l} \in \bdC^{m_{l\mns 1} \tms m_l}$~ is of full
rank since ~$S_{l\mns 1,l} \,\bdz = \bdo$~ with ~$\bdz \ne \bdo$~ will lead to
~$Y_l \,\bdz \in \cN{(A-\la I)^{l\mns 1}} \,=\,\cR{\blb Y_1,\cdots,
Y_{l\mns 1}\brb}$, ~contradicting to the linear independence of ~columns of
~$\blb Y_1,\cdots,Y_l \brb$.
\hfill {\LARGE $\Box$}

The matrix ~$\la_m I + S$~ in (\ref{auus0}) is called a
{\bf local staircase form} of ~$A$~ associated with $\la$. 
~The matrix ~$S$~ is called a
{\bf staircase nilpotent matrix} associated with eigenvalue ~$\la$.
~Writing ~$Y = \blb Y_1,\cdots,Y_k \brb$, ~we call the array
~$(\la, Y, S)$~ as in (\ref{auus0}) \label{'stc3'}
a {\bf staircase eigentriplet}
of ~$A$~ associated with Weyr characteristic ~$\{m_1\geq m_2 \geq \cdots\}$.
~It is called {\bf unitary-staircase eigentriplet} of ~$A$~
if ~$Y$~ is a unitary matrix.
~The unitary-staircase form is often preferable to Jordan Canonical
Form itself since the columns of ~$Y = \blb Y_1, \cdots, Y_k \brb$~
in (\ref{auus0}) form an orthonormal basis for the invariant
subspace of ~$A$~ associated with ~$\la$.

Let ~$\eig{A} = \{\la_1,\cdots,\la_l\}$. 
~Lemma~\ref{lem:scf} lead to the existence of a unitary matrix
~$U\in \bdC^{n\tms n}$~ satisfying \cite{golub-wilkinson,kubl,ruhe-70-bit}
\begin{equation} \label{scdec}
A = U T U^\h, \mbox{\ \ where \ \ } T =
\left[ \mbox{\scriptsize $\begin{array}{cccl}
\la_1 I + S_1 & T_{12}     & \cdots & T_{1l} \\
            & \la_2 I+S_2  & \ddots &  \;\;\; \vdots \\
            &            & \ddots &  T_{l\mns 1,l} \\
            &            &        &  \la_l I + S_l \end{array}$} \right]
\end{equation}
The matrix ~$T$~ in (\ref{scdec}) is called a {\bf staircase
form}\index{staircase form!of a matrix} of ~$A$~ and
the matrix factoring ~$UTU^\h$~ is called a
{\bf unitary-staircase decomposition}
of ~$A$.
~A staircase decomposition of matrix ~$A$~ can be converted to Jordan
decomposition via a series of similarity transformations
\cite{golub-wilkinson,kagstrom-ruhe,kubl,ruhe-70-bit}.

\vspace{-5mm}
\subsection{The notion of the numerical Jordan Canonical Form}
\label{sec:ajcfnotion}
\vspace{-5mm}

\begin{figure}[ht]
\begin{center}
\epsfig{figure=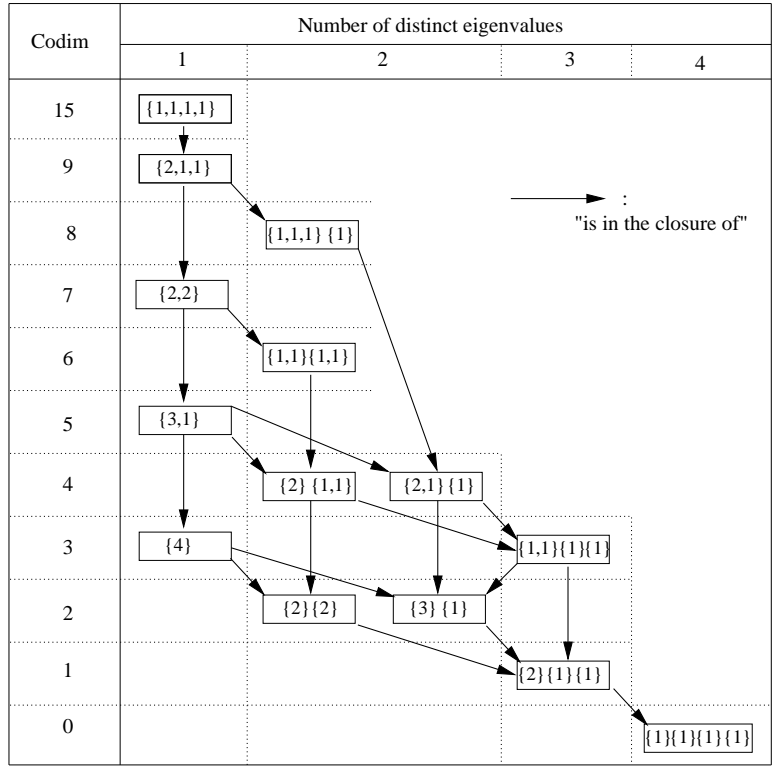,height=4in,width=4in}
\end{center}
\vspace{-3mm} \caption{\footnotesize Stratification of Jordan structures for
~$4\times 4$~ matrices \cite{eek2}. ~General stratification
graphs can be drawn automatically using software package
{\tt StratiGraph} \cite{ejk01,joh03}.} \label{fig:strat}
\end{figure}

Corresponding to a fixed set of ~$k$~ integer partitions
~$\{n_{i1} \ge n_{i2} \ge \cdots\}$~ for ~$i=1,\cdots k$
~with ~$\sum_{i,j} n_{ij} \,=\, n$,
~the collection of all ~$n\times n$ ~matrices with ~$k$~ distinct eigenvalues
associated with Segre characteristics ~$\{n_{ij}\}_{j=1}^\infty$~ for ~$i = 1,
\cdots, k$ ~forms a manifold, known as
a {\bf matrix bundle} originated by A. I. Arnold \cite{arnold}.
~This bundle has a codimension that can be represented
in terms of Segre/Wyre characteristics \cite{arnold,dem-edel}
\begin{equation} \label{bdldim}
\mbox{$\sum_{i=1}^k \Big(-1+\sum_{j=1}^\infty (2j-1)n_{ij}  \Big) ~~\equiv~~
\sum_{i=1}^k \Big(-1 + \sum_{j=1}^\infty m_{ij}^2 \Big)$}.
\end{equation}
where ~$\{m_{ij}\}_{j=1}^\infty$ ~for ~$1 \le i \le k$ ~are corresponding
Weyr characteristics.
~When a matrix ~$A$~ belongs to such a bundle, ~it can also be in the
{\em closure} of many bundles with respect to different Segre characteristics.
~In other words, a matrix with certain Jordan structure can be arbitrarily
close to matrices with other Jordan structures.
~For example, matrix deformations
\begin{equation} \label{matdefm}
\begin{array}{c}
\mbox{\tiny $
\left[ \begin{array}{cccc} \la & 1 && \\ & \la &\eps & \\ && \la & \dl \\
&&& \la \end{array}\right]$}
\\ \mbox{\scriptsize bundle codimention = 3} \end{array}
~~\stackrel{\eps \rightarrow 0~}{\mbox{\LARGE $\longrightarrow$}}~~
\begin{array}{c}
\mbox{\tiny $
\left[ \begin{array}{cccc} \la & 1 && \\ & \la & & \\ && \la & \dl \\
&&& \la \end{array}\right]$}
\\ \mbox{\scriptsize bundle codimention = 7} \end{array}
~~\stackrel{\dl \rightarrow 0~}{\mbox{\LARGE $\longrightarrow$}}~~
\begin{array}{c}
\mbox{\tiny $
\left[ \begin{array}{cccc} \la & 1 && \\ & \la && \\ && \la & \\
&&& \la \end{array}\right]$}
\\ \mbox{\scriptsize bundle codimention = 9} \end{array}
\end{equation}
show that a matrix with Segre characteristic ~$\{2,1,1\}$~ is
arbitrarily close to some matrices with Segre characteristic
~$\{2,2\}$, ~which are arbitrarily near certain matrices with Segre
characteristic ~$\{4\}$.
~Let ~${\cal B}(\cdot)$~ denote the matrix bundle with respect to the Segre
characteristics listed in ~$(\cdot)$ ~and ~$\overline{{\cal B}(\cdot)}$
~denote its closure.
~Then (\ref{matdefm}) suggests that
~${\cal B}(\{2,1,1\}) ~\subset~ \overline{{\cal B}(\{2,2\})}$
~and ~${\cal B}(\{2,2\}) ~\subset~ \overline{{\cal B}(\{4\})}$.
~Extensively studied in e.g. \cite{bhm98,dem-edel,eek1,eek2,ejk03},
these closure relationships form a hierarchy or {\em
stratification} of Jordan structures that can be conveniently
decoded by a covering relationship theorem by Edelman,
Elmroth and K{\aa}gstr\"{o}m \cite[Theorem 2.6]{eek2}.
~As an example, Figure~\ref{fig:strat} lists all the Jordan structures 
and their closure stratification for
~$4\times 4$ ~matrices in different Segre characteristics and
codimensions of matrix bundles.

Let ~$\dist{A,\Pi} \,=\, \inf \big\{\, \|A-B\|_F ~\big|~ B \in \Pi \,\big\}$
present the {\em distance} of a matrix ~$A$ ~to a bundle ~$\Pi$.
~When ~$A \in \bdC^{n\tms n}$ ~is near a bundle ~$\Pi$, ~say
~$\Pi = {\cal B}(\{2,2\})$ ~listed in Figure~\ref{fig:strat}, ~then clearly
~$ \dist{A,\Pi} \ge \dist{A,{\cal B}(\{3,1\})} $ ~since
~$\Pi = {\cal B}(\{2,2\}) \subset \overline{{\cal B}(\{3,1\})}$.
~Indeed, matrix ~$A$ ~is automatically as close, or even closer, to ten other
bundles of {\em lower codimensions} 
below ~$\{2,2\}$ ~following the hierarchy.
~In practical applications and numerical computation,
the given matrix ~$A$ ~comes with imperfect data and/or roundoff error.
~We must assume ~$A = \widehat{A} + E$ ~with a perturbation ~$E$ ~of
small magnitude ~$\|E\|_F$ ~on the original matrix ~$\widehat{A}$.
~The main theme of this article is: ~How to compute the Jordan Canonical
Form of the matrix ~$\widehat{A}$ ~accurately from its inexact data ~$A$.

Suppose matrix ~$\widehat{A}$ ~has an exact nontrivial Jordan Canonical Form and
thus belongs to a bundle ~$\Pi$ ~of codimension ~$d$,
~then there is a lower bound ~$\dl > 0$ ~for the distance from ~$A$  ~to any
other bundle ~$\Pi'$ ~of codimension ~$d' \ge d$.
~When ~$A$ ~is the given data of ~$\widehat{A}$ ~with an imperfect accuracy,
that is, ~$A = \widehat{A}+E$ ~with a perturbation ~$E$, ~then ~$A$ ~generically
resides in the bundle ~${\cal B}(\{1\},\{1\},\cdots,\{1\})$ ~of codimension 0.
~As a result, the Jordan structure of ~$\widehat{A}$ ~is lost in exact
computation on ~$A$ ~for its (exact) Jordan Canonical Form.
~However, the original bundle ~$\Pi$ ~where ~$\widehat{A}$ ~belongs to
has a distinct feature:
it is of the highest codimension among all the bundles passing through
the ~$\eps$-neighborhood of ~$A$, ~as long as ~$\eps$ ~satisfies
~$\dist{A,\Pi} < \eps < \dl$.
~Therefore, to recover the desired Jordan structure of ~$\widehat{A}$ ~from
its empirical data ~$A$, ~we first identify the matrix bundle ~$\Pi$ ~of 
the highest codimension in the neighborhood of ~$A$, ~followed by determining 
the matrix ~$\widetilde{A}$ ~on ~$\Pi$ ~that is closest to ~$A$.
~The {\em numerical Jordan Canonical Form} of ~$A$ ~will then be defined
as the exact Jordan Canonical Form of ~$\widetilde{A}$.
~In summary, the notion of the numerical Jordan Canonical Form is
formulated according to the following three principles:

\vspace{-4mm}
\begin{itemize} \parskip-1mm
\item {\em Backward nearness:} ~The numerical Jordan Canonical Form of
~$A$ ~is the exact Jordan Canonical Form of certain matrix ~$\widetilde{A}$
~within a given distance ~$\eps$,  ~namely ~$\|A-\widetilde{A}\,\|_F < \eps$.
\item {\em Maximum codimension}: ~Among all matrix
bundles having distance less than ~$\eps$ ~of ~$A$, ~matrix
~$\widetilde{A}$ ~lies in the bundle ~$\Pi$ ~with the highest
codimension.
\item {\em Minimum distance}: ~Matrix ~$\widetilde{A}$
~is closest to ~$A$ ~among all matrices in the bundle ~$\Pi$.
\end{itemize}

\vspace{-4mm}
\begin{defn}
~For ~$A \in \bdC^{n\tms n}$ ~and ~$\eps >0$, ~let
~$\Pi \subset \bdC^{n\tms n}$ ~be the matrix bundle such that
\[ \codim{\Pi} ~=~ \max\big\{ \codim{\Pi'} ~\big|~ \dist{A,\Pi'} < \eps\big\},
\]
and ~$\widetilde{A} \in \Pi$ ~satisfying
~$\displaystyle \|A-\widetilde{A}\|_F  =  \min_{B \in \Pi} \|A-B\|_F$
~with {\em (}exact{\em )} Jordan decomposition 
~$\widetilde{A} = X J X^{\mns 1}$.
~Then ~$J$ ~is called the {\bf numerical Jordan Canonical Form} of ~$A$
~within ~$\eps$, ~and ~$XJX^{\mns 1}$ ~is called the
{\bf numerical Jordan decomposition} of ~$A$ ~within ~$\eps$.
\end{defn}
\vspace{-2mm}

{\bf Remark:} ~The same three principles have been successfully applied to
formulate other ill-posed problems with well-posed numerical solutions
such as numerical multiple roots \cite{zeng-05} and
numerical polynomial GCD \cite{zeng-uvgcd,zeng-mvgcd}.
~In this section, we shall attempt to determine the structure of the bundle
~$\Pi$ ~with the highest codimension in the neighborhood of ~$A$.
~The iterative algorithm {\sc EigentripletRefine} developed in
\S \ref{s:itref} is essentially used to find the matrix ~$\widetilde{A}$ 
~in the bundle ~$\Pi$ ~which is nearest to matrix ~$A$.
\hfill {\LARGE $\Box$}

There is an inherent difficulty in computing the Jordan
structure from inexact data and/or using floating point arithmetic.
~If, for instance, a matrix is near several bundles of the same codimension with
almost identical distances, then the structure identification may not be a well
determined problem.
~Therefore, occasional failures \cite{edel-ma-00} for computing the numerical 
Jordan Canonical Form can not be completely eliminated.

\vspace{-5mm}
\section{Regularity of a staircase eigentriplet} \label{s:reg}
\vspace{-5mm}

For an eigenvalue ~$\la$~ of matrix ~$A$~ with a fixed Weyr characteristic,
the components ~$U$~ and ~$S$~ of the staircase eigentriplet
in the staircase decomposition ~$A\,U \,=\, U\,(\la I + S)$~
are not unique.
~We shall impose additional constraints for achieving uniqueness which
is important in establishing the well-posedness of
computing the numerical staircase form.

\vspace{-4mm}
\begin{theorem} \label{unisc}
~~Let ~$A\in \bdC^{n \tms n}$ ~and ~$\la \in \eig{A}$ ~of multiplicity
~$m$
~with nonzero Weyr characteristic ~$m_1\geq \cdots \geq m_k$.
~Then for almost all ~$\bdb_{1}, \cdots, \bdb_{m} \in \bdC^n$~
there is a unitary matrix
~$U = \blb \bdu_1,\cdots,\bdu_m \brb \in \bdC^{n \tms m}$~ and a
staircase nilpotent matrix ~$S \in \bdC^{m\tms m}$~ as in {\em (\ref{sl})}
such that
\begin{eqnarray} \label{auus}
\lefteqn{
\left\{ \begin{array}{rcl} 
 A\,U  -  U(\la I + S) & = & O  \\
\bdu_i^\h \,\bdb_j & = & 0  
~~~~\mbox{for every} ~(i,j) \in \Phi_\la \end{array}
\right.} \\
 ~~\mbox{where} & & 
\label{idxs}
\Phi_\la  ~~\equiv~~ 
\big\{\; (i,j) \,\;\Big|\; \mu_{_{l\mns 1}} < i \leq 
\mu_{_l},
\;\; i < j \leq \mu_{_l}, \;\; l = 1,\cdots,k \; \big\} \\
\label{coweyr}
\mbox{and} & & 
\mu_0 = 0, \;\;\;\; \mu_j = m_1 + \cdots + m_j, \;\;\; j = 1, \cdots, k.
\end{eqnarray}
Moreover, 
if there is another unitary matrix ~$\hat{U} = [\hat{\bdu}_1,\cdots,
\hat{\bdu}_m]$~ and a staircase nilpotent matrix ~$\hS$~
that can substitute ~$U$ ~and ~$S$~ in {\em (\ref{auus})}, 
then ~$\hS = \hat{U}^\h (A-\la I) \hat{U}$ ~where
~$\hat{U} = U D$ ~for a diagonal matrix ~$D = \diag{\al_1,\cdots,\al_m}$ 
~with ~$|\al_1|=\cdots =|\al_m|=1$.
\end{theorem}
\vspace{-4mm}

\begin{wrapfigure}[8]{r}{3.0in}
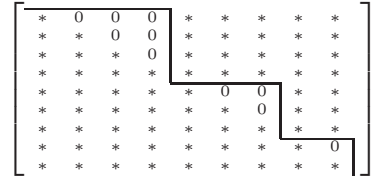

\vspace{-4mm}
\[
\left[ \mbox{\tiny $
\begin{array}{ccccccccc} \cline{1-4}
* & 0 & 0 & 0 & \multicolumn{1}{|c}{*} & * & * & * & * \\
* & * & 0 & 0 & \multicolumn{1}{|c}{*} & * & * & * & * \\
* & * & * & 0 & \multicolumn{1}{|c}{*} & * & * & * & * \\
* & * & * & * & \multicolumn{1}{|c}{*} & * & * & * & * \\ \cline{5-7}
* & * & * & * & * & 0 & 0 & \multicolumn{1}{|c}{*} & * \\
* & * & * & * & * & * & 0 & \multicolumn{1}{|c}{*} & * \\
* & * & * & * & * & * & * & \multicolumn{1}{|c}{*} & * \\ \cline{8-9}
* & * & * & * & * & * & * & * & \multicolumn{1}{c|}{0} \\
* & * & * & * & * & * & * & * & \multicolumn{1}{c|}{*} \\
\end{array}$ } \right]
\]
\vspace{-4mm}
\caption{\footnotesize Index set ~$\Phi_\la$ ~for Weyr characteristic
~$\{4,3,2\}$. ~Every ~$(i,j) \in \Phi_\la$ ~entry is zero.} 
\label{phi}
\end{wrapfigure}
$\mathbf{Remark:}$ 
~The second equation in (\ref{auus}) along with the index set 
~$\Phi_\la$ ~in (\ref{idxs}) means that, for ~$j=1,\cdots,k$, 
the matrix ~$U_j B_j$ ~is lower triangular where 
~$U_j = \blb \bdu_{\mu_{j\mns 1} \pls 1}, \cdots, \bdu_{\mu_{j}} \brb$
~and
~$B_j = \blb \bdb_{\mu_{j\mns 1} \pls 1}, \cdots, \bdb_{\mu_{j}} \brb$.
For example: Let ~$\la$ ~be an eigenvalue 
with Weyr characteristic ~$\{4,3,2\}$, ~we have multiplicity 
~$m = 9$, ~$[\mu_0, \mu_1, \mu_2,
\mu_3] = [0, 4, 7, 9]$, and 
~$\Phi_\la = \{$
{\scriptsize $(1,2)$, $(1,3)$, $(1,4)$, $(2,3)$, $(2,4)$, $(3,4)$,
$(5,6)$, $(5,7)$, $(6,7)$, $(8,9)$} $\}$. 
~The matrix ~$U^\h\, [\bdb_1,\cdots,\bdb_9]$ ~has zero at every 
~$(i,j) \in \Phi_\la$ ~entry. 
~As shown in Figure~\ref{phi}, those zeros are under the
staircase and above the diagonal. \qed

%\begin{proof}
{\bf Proof of Theorem \ref{unisc}}.
~For ~$j=1,\cdots,k$, ~the subspace 
~$\cN{(A-\la I)^j}\cap \cN{(A-\la I)^{j\mns 1}}^\perp$ ~is of dimension 
~$m_j$.
~For almost all vectors ~$\bdb_i$, ~$i=\mu_{j\mns 1} \lpls 1, \cdots, 
\mu_j$, ~the subspace 
\[
\cN{(A-\la I)^j}\cap \cN{(A-\la I)^{j\mns 1}}^\perp \cap
\spn{\bdb_{\mu_{j\mns 1} \pls 2}, \cdots \bdb_{\mu_j}}^\perp \]
is of dimension one and spanned by a unit vector 
~$\bdu_{\mu_{j\mns 1} \pls 1}$ ~which is unique up to a unit constant
multiple.
~After obtaining ~$\bdu_{\mu_{j\mns 1} \pls 1}, \cdots 
\bdu_{\mu_{j\mns 1} \pls l}$, ~the subspace 
\[
\cN{(A-\la I)^j}\cap \cN{(A-\la I)^{j\mns 1}}^\perp \cap
\spn{\bdu_{\mu_{j\mns 1} \pls 1}, \cdots 
\bdu_{\mu_{j\mns 1} \pls l},\bdb_{\mu_{j\mns 1} \pls l \pls 2}, 
\cdots \bdb_{\mu_j}}^\perp \]
is of dimension one and spanned by ~$\bdu_{\mu_{j\mns 1} \pls l \pls 1}$
~which is again unique up to a unit constant multiple.
~Therefore, we have a unitary matrix 
~$U_j = \blb \bdu_{\mu_{j\mns 1} \pls 1}, \cdots, \bdu_{\mu_{j}} \brb$, 
whose columns satisfy the second equation in (\ref{auus}) and span the subspace
~$\cN{(A-\la I)^j}\cap \cN{(A-\la I)^{j\mns 1}}^\perp$ ~for 
~$j=1,\cdots,k$. 
~These unitary matrices uniquely determines ~$S_{ij} = U_i^\h (A-\la I) U_j$ 
~in (\ref{sl}). 
~It is straightforward to verify (\ref{auus}) 
for ~$U = \blb U_1, \cdots, U_k \brb = \blb \bdu_1,\cdots, \bdu_{m} \brb$.
\qed

One of the main components of our algorithm is an iterative refinement of 
the eigentriplet ~$(\la,U,S)$ ~using the Gauss-Newton iteration.
~For this purpose we need to construct a system of analytic equations for the 
eigentriplet ~$(\la, U, S)$~ specified in Theorem \ref{unisc}.
~If the matrix ~$A$ ~and the eigenvalue ~$\la$ ~are real, it is 
straightforward to set up the system using equations in ~$(\ref{auus}$ ~along
with the orthogonality equaiton ~$U^\top U-I=O$.
~When either matrix ~$A$~ or eigenvalue ~$\la$~ is complex, however,
the unitary constraint ~$U^\h U - I = O$ ~is not analytic. 
~One way to circumvent this difficulty is converting (\ref{auus}) and
~$U^\h U - I = O$ ~to real equations by splitting ~$A$ ~and the eigentriplet
~$(\la, U, S)$ ~into real and imaginary parts.
%
%~The system of real equations 
%\[ \left\{ \begin{array}{rcl}
%\hat{A}\hat{U}-\check{A}\check{U} - \hat{U}(\hat{\la} I+\hat{S}) + 
%\check{U}(\check{\la} I +\check{S}) & = & O \\
%\hat{A}\check{U}+\check{A}\hat{U} - \hat{U}(\check{\la} I+\check{S}) - 
%\check{U}(\hat{\la} I +\hat{S}) & = & O \\
%\hat{\bdu}_i^\top \hat{\bdb}_j + \check{\bdu}_i^\top\check{\bdb}_j 
%& = & \bdo \\
%\hat{\bdu}_i^\top \check{\bdb}_j - \check{\bdu}_i^\top \hat{\bdb}_j 
%& = & \bdo, ~~~~(i,j) \in \Phi_\la \\
%\hat{U}^\top \hat{V} + \check{U}^\top\check{V}
%& = & I  \\
%\hat{U}^\top \check{V} - \check{U}^\top \hat{V}
%& = & O 
%
%\end{array} \right.
%\]
%where ~$(\hat{\cdot})$ ~and ~$(\check{\cdot})$
%~represents the real and imaginary parts of ~$(\cdot)$ ~respectively.
%
~The resulting system of real equations would be real analytic.

Alternatively, we developed a simple and effective strategy to overcome this 
difficulty by a two step approach. 
~As an initial approximation, a staircase eigentriplet ~$(\la, Y, S)$ ~is 
computed and ~$Y$ ~does not need to be unitary. 
~We replace the unitary constraint ~$U^\h U - I$~ with a nonsingularity 
requirement
\begin{equation} \label{cybde}
C^\h \,Y ~~=~~ \left[ \mbox{\tiny
$\begin{array}{cccc} 
1 & 0 & \cdots & 0 \\
* & 1 & \ddots & \vdots \\ \vdots & \ddots & \ddots & 0 \\
* & \cdots & * & 1 \end{array}$} \right]_{m\tms m}
%\, \blb \bdc_1, \bdc_2, \cdots, \bdc_i \brb^\h \,\bdy_i ~~=~~
%[\underbrace{0,0,\cdots,0}_{i-1},1]^\h
%\mbox{\ \ for \ } 1 \leq i \leq m,
\end{equation}
via a constant matrix
~$C = \{\bdc_1,\cdots,\bdc_m\} \subset \bdC^n$.
~The solution ~$Y = [\bdy_1,\cdots,\bdy_m]$~ to the equation (\ref{cybde})
combined with (\ref{auus}) is a nonsingular matrix whose columns span the
invariant subspace of ~$A$ ~associated with ~$\la$.
~Then, at the second step, ~$Y$ ~can be stably orthogonalized to 
~$U = \blb \bdu_1,\cdots,\bdu_m\brb$~ and provide the solution of (\ref{auus}).
~If necessary, we repeat the process as refinement by replacing 
~$C$ ~with ~$U$ ~in ~(\ref{cybde}) along with (\ref{auus}) and solve for 
~$Y$ ~again using the previous eigentriplet results as the initial iterate.
~In the spirit of Kahan's well-regarded ``twice is enough'' observation%
~\cite[p. 110]{ParlettBook}, this reorthogonalization never needs the 
third run.

In general, let ~$\{\bdb_1,\cdots,\bdb_m\} \subset \bdC^n$~ be the set of
predetermined random complex vectors as in Theorem~\ref{unisc}.
~A second set of complex vectors ~$\bdc_1,\cdots,\bdc_m \in \bdC^n$~ will
also be chosen to set up the overdetermined
quadratic system
\begin{equation}\label{ayys}
\left\{ \begin{array}{rcll}
(A - \la I) Y  & = &  Y\,S  & \\
\blb \bdc_1,\cdots,\bdc_i \brb^\h\, \bdy_i & = & [0,\cdots,0,1]^\h,
& \mbox{\ \ for \ } 1 \leq i \leq m \\
 \bdb_j^\h\,\bdy_i & = & 0, &\mbox{\ \ for \ \ } (i,j) \in \Phi_\la
\end{array} \right.
\end{equation}
where ~$\Phi_\la$ ~is defined in  (\ref{idxs}).
~There are ~$\eta$~ equations and ~$\zeta$~ unknowns in (\ref{ayys}) where
\begin{equation} \label{sysize}
\mbox{$ \eta ~~=~~ nm + \frac{m^2}{2} + \frac{1}{2}\sum_{j=1}^k m_j^2
\mbox{\ \ \ and \ \ }
\zeta ~~=~~ 1 + nm + \sum_{i<j} m_i m_j$}
\end{equation}
with a difference ~$\eta-\zeta = -1 + \sum m_j^2$.
~Let
\begin{equation} \label{cF}
\bdf(\la, Y, S) ~~=~~ \mbox{\footnotesize $\left[ \begin{array}{c}
\big( (A - \la I) Y - Y S \big)\,\bde_1 \\
\;\;\;\;\;\; \vdots \\
\big( (A - \la I) Y - Y S \big) \, \bde_m \\
\llbracket \bdc_j^\h \bdy_i -\dl_{ij}\rrbracket \\
\llbracket \bdb_{j}^\h \bdy_i \rrbracket
\end{array} \right]$}
\end{equation}
where ~$\dl_{ij}$~ is the Kronecker delta,
~$\llbracket \bdc_j^\h \bdy_i - \dl_{ij}\rrbracket$~ 
and ~$\llbracket \bdb_{j}^\h \bdy_i \rrbracket$~
denote vectors of components
~$\big\{\, \bdc_{j}^\h \bdy_i - \dl_{ij} \, \big|\, 1 \le i \le m, \;
j \le i \big\}$~ and ~$\big\{\, \bdb_{j}^\h \bdy_i\, \big|\,
(i,j) \in \Phi_\la \,\big\}$ ~respectively,  ordered by the rule where
~$(i,j)$~ precedes ~$(i',j')$~ if ~$i < i'$, ~or ~$i=i'$~ with ~$j < j'$.
~The ~$\zeta$~ unknowns in eigentriplet ~$(\la, Y, S)$~ are ordered in a
vector form
\begin{equation} \label{triporder}
 (\la, \bdy_1^\top, \cdots, \bdy_m^\top, \bds^\top )^\top \end{equation}
where ~$\bds$~ is the column vector consists of the entries of ~$S$~
in the order illustrated in the following example for the Weyr characteristic
~$\{3 \geq 2 \geq 1\}$:
\begin{eqnarray}  \label{S321}
S & = & \left[ \mbox{\scriptsize
$\begin{array}{cccccc} \cline{1-3}
0& 0& 0& \multicolumn{1}{|c}{s_{14}} & s_{15} & s_{16}  \\
0& 0& 0& \multicolumn{1}{|c}{s_{24}} & s_{25} & s_{26}  \\
0& 0& 0& \multicolumn{1}{|c}{s_{34}} & s_{35} & s_{36}  \\ \cline{4-5}
0& 0& 0&        0&        0& \multicolumn{1}{|c}{s_{46}}  \\
0& 0& 0&        0&        0& \multicolumn{1}{|c}{s_{56}}  \\ \cline{6-6}
0& 0& 0&        0&        0& \multicolumn{1}{c|}{0} \end{array}$} \right]
\\ \nonumber 
\bds^\top  & = &[s_{14}, s_{24}, s_{34},\;\;\; s_{15}, s_{25}, s_{35},\;\;\;
s_{16}, s_{26}, s_{36}, s_{46},s_{56}] 
\end{eqnarray}
With this arrangement, the Jacobian ~$J(\la,Y,S)$~ of ~$\bdf(\la, Y, S)$~ is
an ~$\eta \times \zeta$~ matrix. 

\vspace{-4mm}
\begin{theorem} \label{t:cplxreg}
~~Let ~$\la$~ be an ~$m$-fold eigenvalue of ~$A \in \bdC^{n\tms n}$~
associated with nonzero Weyr characteristic
~$\{m_1 \geq \cdots \geq m_k\}$.
~Then for almost all vectors ~$\bdb_1,\cdots,\bdb_m,
\bdc_1,\cdots,\bdc_m \in \bdC^n$, ~there is a unique pair of matrices
~$Y = [\bdy_1,\cdots,\bdy_m]\in \bdC^{n\tms m}$~ and  ~$S \in \bdC^{m \tms m}$
~where ~$S$~ is a staircase nilpotent
matrix in the form of \,{\em (\ref{sl})}\, such that
the staircase eigentriplet ~$(\la,Y,S)$~ satisfies the system
{\em (\ref{ayys})}.
~Moreover, the Jacobian ~$J(\cdot,\cdot,\cdot)$~ of
~$\bdf(\cdot,\cdot,\cdot)$~ in {\em (\ref{cF})} is of full column rank
at ~$(\la, Y, S)$.
\end{theorem}
\vspace{-4mm}

\prf 
The subspace ~$\cN{(A-\la I)^j}$ ~is of
dimension ~$\mu_j$ ~for ~$j=1,\cdots,k$. 
~For each ~$l \in \{\mu_{j\mns 1}+ 1, \cdots, \mu_j\}$, ~the subspace 
~$\cN{(A-\la I)^j} \bigcap \spn{\bdc_1,\cdots,\bdc_{l\mns 1},
\bdb_{l\pls 1},\cdots, \bdb_{\mu_j}}^\perp$ 
~is of dimension one and spanned by the unique vector ~$\bdy_l$ ~with
~$\bdc_l^\h \bdy_l = 1$.
~Theirfore vectors ~$\bdy_1,\cdots,\bdy_m$ ~are uniquely defined so that
~$\cR{\blb Y_1,\cdots, Y_j\brb} = \cN{(A-\la I)^j}$ ~for ~$j = 1, \cdots, k$
~where 
~$Y_i = \blb \bdy_{\mu_{i\mns 1}\pls 1}, \cdots \bdy_{\mu_i} \brb$, 
~$i=1,\cdots, k$.
~Moreover, ~$\bdy_i^\h\,\blb \bdc_1,\cdots,\bdc_i] = \blb 0, \cdots, 0, 1\brb$
~and ~$\bdb_j^\h \bdy_i = 0$ ~for ~$(i,j) \in \Phi_\la$.
~It is straightforward to verify 
~$A\blb Y_1,\cdots, Y_k \brb = \blb Y_1,\cdots, Y_k \brb (\la I + S)$
~for ~$S$ ~being a nilpotent staircase matrix in the form of (\ref{sl})
with uniquely determined blocks
\[ \mbox{\scriptsize $
\left[ \begin{array}{l} S_{1j} \\ ~~\vdots \\ S_{j\mns 1, j} \end{array} 
\right]$}
~=~ \blb Y_1,\cdots, Y_{j\mns 1} \,\brb^+ (A-\la I)\, Y_j
~~~\mbox{for} ~j=2,\cdots,k. \]
Consequently, the matrix pair ~$(Y,S)$ ~satisfying (\ref{ayys}) exists
and is unique for almost all ~$\bdc_1,\cdots,\bdc_m,\bdb_1,\cdots,\bdb_m
\in \bdC^n$.

We now prove the Jacobian ~$J(\la,Y,S)$~ of ~$\bdf(\la,Y,S)$~ in
(\ref{cF}) is of full rank at a staircase eigentriplet ~$(\la, Y, S)$.
~The Jacobian ~$J(\la,Y,S)$~ can
be considered a linear transformation which maps ~$(\sg, Z, T)$~ into
~$\bdC^\zeta$, ~where ~$\sg \in \bdC$, ~$Z \in \bdC^{n \tms m}$~ and
~$T$~ is a nilpotent staircase matrix of ~$m\times m$~
relative to
the nonzero Weyr characteristic ~$\{m_1 \geq \cdots \geq m_k\}$.
~We partition ~$T$~ with blocks ~$T_{ij} \in \bdC^{m_i \tms m_j}$~
in the same way as we partition ~$S = S_k$~ in (\ref{sl}) for ~$l=k$.
~Assume ~$J(\la,Y,S)$~ is rank-deficient.
~Then there is a triplet ~$(\sg, Z, T) \ne (0, O, O)$~ such that
~$J(\la,Y,S) [\sg, Z, T] = \bdo$, namely
\begin{eqnarray}
(A - \la I) Z & = & \sg Y + ZS + Y T
\label{azyzsyt} \\
\blb \bdc_1, \cdots, \bdc_i \brb^\h \bdz_i & = & \bdo, \;\;\;
i = 1, \cdots, m, \label{cizicjzi} \\
\bdb_j^\h \bdz_i & = & 0,  \;\;\; (i,j) \in \Phi_\la.
\label{bjzi}
\end{eqnarray}
Using
~$Y_j = \blb \bdy_{\mu_{_{j\mns 1}}\pls 1}, \cdots, \bdy_{\mu_{_j}} \brb$ ~and
~$Z_j = \blb \bdz_{\mu_{_{j\mns 1}}\pls 1}, \cdots, \bdz_{\mu_{_j}} \brb$ ~for
~$j=1,\cdots,k$, ~we have
\begin{equation} \label{alv}
\left\{ \begin{array}{rcl}
(A-\la I) Z_1 & = & \sg Y_1  \\
(A-\la I) Z_i & = & \sg Y_i +
\sum_{j=1}^{i\mns 1} (Z_j S_{ji} + Y_j T_{ji}), ~~~~i = 2, \cdots k.
\end{array} \right.
\end{equation}
from (\ref{azyzsyt}).
~A simple induction using (\ref{alv}) 
leads to ~$(A-\la I)^{j\pls 1} Z_j = O$, ~for ~$j = 1, \cdots, k$.
~Namely, vectors ~$\bdz_1, \cdots, \bdz_m$~ all belong to the invariant
subspace of ~$A$~ associated with ~$\la$, ~and thus ~$Z = YE$ ~holds for
certain ~$E \in \bdC^{m\tms m}$.

Also by a straightforward induction we have
\begin{equation} \label{alx}
(A - \la I)^l Y_{l\pls 1} ~~=~~ Y_1 S_{12}S_{23}\cdots S_{l,l\pls 1},
~~~\mbox{for}~~ l = 1, \cdots, k-1.
\end{equation}
We claim that
\begin{equation} \label{aly}
(A-\la I)^l Z_l ~~=~~ l \sg Y_1 S_{12}S_{23}\cdots S_{l\mns 1,l},
\mbox{\ \ for each \ \ } l = 1, \cdots, k.
\end{equation}
This is true for ~$l=1$~ because of (\ref{alv}).
~Assume (\ref{aly}) is true for ~$l\leq j-1$.
~Then by (\ref{alx})
\begin{eqnarray*}
(A - \la I)^{j} Z_{j} & = & (A - \la I)^{j\mns 1} (A - \la I) Z_j
~ = ~ (A - \la I)^{j\mns 1} \big[\sg Y_j +
\mbox{$\sum_{i=1}^{j\mns 1}$} (Z_i S_{ij} + Y_i T_{ij}) \big] \\
& = & \sg (A - \la I)^{j\mns 1} Y_j +
(A-\la I)^{j\mns 1} Z_{j\mns 1} S_{j\mns 1,j} 
%& = & \sg Y_1 \prod_{i=1}^{j\mns 1} S_{i,i\pls 1} +
%(j-1) \sg Y_1 \prod_{i=1}^{j\mns 1} S_{i,i\pls 1} +
%& & \;\;\;\;\;\;\;\;\; +
%\sum_{i=1}^{j\mns 2} (A-\la I)^{j\mns 1\mns i}
%[i\sg Y_1 S_{12}\cdots S_{i\mns 1,i}] S_{ij} \\
~ = ~ j \sg Y_1 S_{12} S_{23}\cdots  S_{j\mns 1,j}
\end{eqnarray*}
since, again, ~$(A-\la I)^{j\mns 1} Y_i = O$ ~for ~$i \le j-1$.
~Thus (\ref{aly}) holds for ~$l=1,\cdots,k$.

Since, ~$(A-\la I)^k Y \,=\, O$,
~hence
~$(A-\la I)^k Z = (A-\la I)^k Y E \,=\, O$
~from~ ~$Z = YE$.
~From (\ref{aly}), we have 
~$(A-\la I)^k Z_k \,=\, k \sg Y_1 S_{12} S_{23}\cdots S_{k\mns 1,k} = O$,
~By Lemma \ref{lem:scf}, ~$S_{12} S_{23}\cdots S_{k\mns 1,k}$~ is of full
rank.
~Consequently, ~$(A-\la I)^l Z_l = O$ ~by (\ref{aly}), namely
~$\cR{Z_l} \subset \cN{(A-\la I)^l}$ ~for ~$l=1,\cdots,k$.
~Therefore, for every ~$i \in \{\mu_{i\mns 1}+1,\cdots,\mu_i\}$, ~$\bdz_i$
~is in 
\[ \cN{(A-\la I)^l} \bigcap 
\spn{\bdc_1,\cdots,\bdc_i,\bdb_{i\pls 1},\cdots,\bdb_{\mu_i}}^\perp 
~=~ \{\bdo\}.
\]
for ~$i = 1, \cdots, k$, ~implying ~$Z = O$. 
~The equation (\ref{azyzsyt}) them implies ~$YT= O$~ and
thus ~$T=O$~ since ~$Y$~ is of full column rank.
~Consequently, ~$J(\la,Y,S)$~ is of full column rank. \hfill {\LARGE $\Box$}

The component ~$Y$~ in the staircase eigentriplet ~$(\la, Y, S)$~
satisfying (\ref{ayys}) is a unitary matrix for a particular 
~$\blb \bdc_1,\cdots,\bdc_m \brb = Y$.
~This will be achieved in our eigentriplet refinement process.

\vspace{-5mm}
\section{The numerical staircase eigentriplet and its sensitivity}
\label{sec:sens}
\vspace{-5mm}

Consider an $n \times n$ ~complex matrix ~$A$ ~along with a fixed partition
~$\{\,m_1 \ge m_2 \ge \cdots \ge m_k > 0\,\}$~ of integer (multiplicity)
~$m > 0$.
~Let the vector function ~$\bdf(\la,Y,S)$ ~be defined defined in (\ref{cF})
with respect to fixed vectors ~$\bdb_1,\cdots,\bdb_m$~ and
~$\blb \bdc_1,\cdots,\bdc_m \brb = \blb \bdy_1,\cdots,\bdy_m \brb$ ~and
~$J(\cdot,\cdot,\cdot)$ ~is its Jacobian.
~An array ~$(\la,Y,S) \in \bdC \times \bdC^{n\tms m} \times \bdC^{m\tms m}$~
is called a {\bf numerical unitary-staircase eigentriplet} of ~$A$~
with respect to nonzero Weyr characteristic
~$\{\,m_1 \ge m_2 \ge \cdots \ge m_k \,\}$~
if ~$(\la,Y,S)$~ satisfies 
~$J(\la,Y,S)^\h \bdf(\la,Y,S) = \bdo$, ~a necessary condition for
~$\big\|\bdf(\cdot,\cdot,\cdot)\big\|_2$ ~to reach a local minimum
at ~$(\la,Y,S)$.
~The requirement 
~$\blb \bdc_1,\cdots,\bdc_m \brb = \blb \bdy_1,\cdots,\bdy_m \brb$ 
can be satisfied in our refinement algorithm that is to be
elaborated in \S \ref{s:itref}.

If ~$A$~ possesses a numerical unitary-staircase eigentriplet ~$(\la, Y, S)$~
with a small residual
\begin{equation} \label{resdef}
  \rho \,~=~\,\|\,A\,Y\,-\,Y\,(\la I + S)\,\|_F\big/\|A\|_F,
\end{equation}
then, letting ~$Z$~ be a unitary complement of ~$Y$, ~it is
straightforward to verify that
\begin{equation} \label{nearA}
\hat{A} \;~=~\; [\mbox{\scriptsize $Y,\,Z$}]
\left[\mbox{\scriptsize $\begin{array}{cc} \la I + S &
Y^\h A Z \\ O & Z^\h A Z \end{array}$} \right]
\left[ \mbox{\scriptsize $\begin{array}{c} Y^\h \\ Z^\h \end{array}$} \right]
~~=~~ Y (\la I + S) Y^\h + Y Y^\h A Z Z^\h + Z Z^\h A Z Z^\h
\end{equation}
possesses ~$(\la,Y,S)$~  as its {\em  exact} unitary-staircase eigentriplet
and the distance
\begin{eqnarray*}
\big\|\, A - \hat{A} \,\big\|_F & = & \big\|\, (A - \hat{A})[Y,Z] \,\big\|_F 
~~=~~ \big\|\, (A-\hat{A}) Y \,\|_F + \big\|\,(A-\hat{A})Z \,\big\|_F \\
& = & \big\|\, A Y - Y (\la I + S) \,\big\|_F + 
\big\|\, AZ - (YY^\h A Z + ZZ^\h AZ) \,\big\|_F \\ 
& = & \big\|\, A Y - Y (\la I + S) \,\big\|_F ~~=~~ \rho \,\|A\|_F 
\label{resdis}
\end{eqnarray*}
is small. 
~We now derive the well-posedness and the sensitivity measurement in a heuristic
manner.
~Let ~$(\la,Y,S)$~ be an ~$m$-fold numerical
unitary-staircase eigentriplet of ~$A$~ with residual
~$\bdq = \bdf(\la,Y,S)$~ for ~$\bdf$~ defined in (\ref{cF}) via certain
auxiliary vectors ~$\bdb_1,\cdots,\bdb_m$~ and ~$\bdc_1,\cdots,\bdc_m$.
~To analyze the effect of perturbation on matrix ~$A$, ~let
~$\bdg(A,\la,Y,S)$~ denote the same vector function ~$\bdf(\la,Y,S)$~ in
(\ref{cF}) where ~$A$~ is now considered as a variable.
~When ~$A$~ becomes ~$\tilde{A} = A + E$~ by adding a matrix
~$E$~ of small norm, denote
~$(\tilde{\la}, \tilde{Y}, \tilde{S})$~ as a numerical unitary eigentriplet
of ~$\tilde{A}$.
~Let us estimate the asymptotic bound of error
\[
\big\|\llbracket \la,Y,S\rrbracket -
\llbracket \tilde{\la}, \tilde{Y}, \tilde{S}\rrbracket \big\|_2 \,~\equiv~\,
\sqrt{|\la-\tilde{\la}|^2+\|Y-\tilde{Y}\|_F^2 + \|S-\tilde{S}\|_F^2}
\]
where ~$\llbracket \la,Y,S \rrbracket$~ and 
~$\llbracket \tilde{\la}, \tilde{Y}, \tilde{S}\rrbracket$~ denote the
vector forms of ~$(\la,Y,S)$~ and ~$(\tilde{\la}, \tilde{Y}, \tilde{S})$~
respectively according to the rule given in (\ref{triporder}).
~Write ~$\bdg(A,\la,Y,S) = \bdq$.
~Since ~$\|\bdg(\tilde{A},\tilde{\la},\tilde{Y},\tilde{S})\|_2$~ is
the local minimum in a neighborhood of ~$(\tilde{\la},\tilde{Y},\tilde{S})$,
~we have
\[ \|\bdg(\tilde{A},\tilde{\la},\tilde{Y},\tilde{S})\|_2 \,~\le~ \,
\|\bdg(\tilde{A},\la,Y,S)\|_2 \, ~\le~ \, \|E\,Y\|_F  + \|\bdq\|_2 \,~\le~\,
\|E\|_F + \|\bdq\|_2 \]
for small ~$\|E\|_F$~ and ~$\|\bdq\|_2$.
~Moreover,
\[ \|\bdg(A,\tilde{\la},\tilde{Y},\tilde{S})\|_2
\,~\le~\, \|\bdg(\tilde{A},\tilde{\la},\tilde{Y},\tilde{S})\|_2 + 
\|E\,\tilde{Y}\|_F \,~\le~\, 2 \|E\|_F + \|\bdq\|_2. \]
In other words,
\begin{eqnarray*}
\|J(\la,Y,S)(\llbracket \la,Y,S\rrbracket- 
\llbracket \tilde{\la},\tilde{Y},\tilde{S}\rrbracket)\|_2
& = &  \|\bdf(\la,Y,S) - \bdf(\tilde{\la},\tilde{Y},\tilde{S})\|_2
\;+\; h.o.t. \\
& \le & 2 \|E\|_F \;+\; 2 \|\bdq\|_2 \;+\;  h.o.t. \end{eqnarray*}
where ~$J(\cdot,\cdot,\cdot)$~ is the Jacobian of ~$\bdf(\cdot,\cdot,\cdot)$~
and ~$h.o.t$~ represents the higher order terms of ~$\|E\|_F+\|\bdq\|_2$.
~Let ~$\sg_{\min}(\cdot)$~ be the smallest singular value  \label{'sgm'}
of matrix ~$(\cdot)$.
~Then ~$\sg_{\min}\big(J(\la,Y,S)\big)$~ is strictly positive by
Theorem \ref{t:cplxreg} and
\begin{eqnarray*} \sg_{\min}(J(\la,Y,S))\,
\bnorm{\llbracket \la,Y,S\rrbracket- 
\llbracket \tilde{\la},\tilde{Y},\tilde{S}\rrbracket}  &\le &
\bnorm{J(\la,Y,S)(\llbracket \la,Y,S\rrbracket - 
\llbracket \tilde{\la},\tilde{Y},\tilde{S}\rrbracket)} \\
& \le &  2 \|E\|_F + 2 \|\bdq\|_2 + h.o.t.
\end{eqnarray*}
where ~$h.o.t$~ represents the higher order terms of ~$\|E\|_F$.
~This provides an asymptotic bound
\begin{equation}
|\,\la - \tilde\la\,| ~\;\le~~
\bnorm{\llbracket \la,Y,S\rrbracket - 
\llbracket \tilde{\la},\tilde{Y},\tilde{S}\rrbracket}
 ~~\le~~  \frac{2}{\sg_{\min}\big(J(\la,Y,S)\big)}
\left( \|E\|_F +  \|\bdq\|_2 \right),
\label{tripbound}
\end{equation}
and the finite positive real number
\begin{equation} \label{tripcond}
\kappa (\la,Y,S) ~\equiv~ 2\,\sg_{\min}\big(J(\la,Y,S)\big)^{\mns 1} ~=~
2\|J(\la,Y,S)^+\|_2
\end{equation}
serves as a condition number of the unitary-staircase eigentriplet that
measures its sensitivity with respect to perturbations on matrix ~$A$.

\vspace{-3mm}
\begin{defn} \label{'stcond'}
~~Let ~$(\la,Y,S)$~ be a numerical unitary-staircase eigentriplet
of ~$A \in \bdC^{n \tms n}$~ as a regular orthogonal solution to
the system ~$\bdf(\cdot,\cdot,\cdot)=\bdo$~ corresponding to auxiliary
vectors ~$\bdb_1,\cdots,\bdb_m$~ and ~$\blb \bdc_1, \cdots,\bdc_m
\brb = Y$~ in {\em (\ref{cF})}.
~Let ~$J(\cdot,\cdot,\cdot)$~ be the Jacobian of ~$\bdf(\cdot,\cdot,\cdot)$.
~Then we call ~$\kappa(\la,Y,S) \equiv 2\|J(\la,Y,S)^+\|_2$~ the
{\bf staircase condition number} for the eigentriplet.
\end{defn}
\vspace{-3mm}

{\bf Remark.} ~~The arithmetic mean of
an eigenvalue cluster is often used as an approximation to a multiple 
eigenvalue.
~Let ~$\la$~ be an ~$m$-fold eigenvalue of ~$A$~ with an orthonormal basis 
matrix ~$Y$~ for the invariant subspace.
~The perturbed matrix ~$A+E$~ has a cluster of eigenvalues around ~$\la$.
~Chatelin \cite[pp.155--156]{chatelin} established the bound on the arithmetic
mean ~$\hat{\la}$~ as
\begin{equation} \label{clcond}
 \big| \hat{\la} - \la \big| ~~\le~~ \big\|(X^\h Y)^{\mns 1} \big\|_2 \|E\|_2
\end{equation}
for small ~$\|E\|_2$, ~where ~$X$~ is a matrix whose columns form a basis
for the invariant subspace of ~$A^\h$.
~We call ~$\big\|(X^\h Y)^{\mns 1}\big\|_2$~ the 
{\bf cluster condition number} of ~$\la$.
~From our computing experiments, the cluster condition number can be
substantially
larger than the staircase condition number as shown in the following example.
\hfill {\LARGE $\Box$}

\vspace{-3mm}
\begin{example}  \label{e:sens} \em
~Matrix%
%\footnote{The Matlab code of our algorithm and test scripts for all the 
%numerical examples can be 
%accessed at {\tt http://www.neiu.edu/$\sim$zzeng/numjcf.htm}}
\[ A \,~=~\,
\mbox{\tiny $
\left[ \begin{array}{rrrrrrrrrrrrrrrrrrrr}
 1&  0& 20& 93&  0& 71& 34&  6&\mbox{-}20&  3& 31&\mbox{-}14&  0&  0&\mbox{-}19& 14&
\mbox{-}11&  0&  3& \mbox{-}6 \\
17&  7&\mbox{-}32&\mbox{-}84&  3&\mbox{-}69&\mbox{-}33& \mbox{-}9& 21& \mbox{-}3
&\mbox{-}30& 17& \mbox{-}2& \mbox{-}4& 15&\mbox{-}12&  9&  0& \mbox{-}3&  9 \\
10&  5& 40&247&  0&193& 92& 17&\mbox{-}58&  9& 83&\mbox{-}21& \mbox{-}5&  0&
\mbox{-}48
& 27&\mbox{-}25&  0&  9&\mbox{-}17 \\
\mbox{-}7& \mbox{-}3&  0&\mbox{-}39& \mbox{-}3&\mbox{-}34&\mbox{-}19& \mbox{-}4& 13
&  0&\mbox{-}15& \mbox{-}1&  0&  0&  7&  1&  2&  0&  0&  4 \\
\mbox{-}6&  0& 62&307& \mbox{-}1&248&118& 26&\mbox{-}77& 12&106&\mbox{-}30& 
\mbox{-}3
&  4&\mbox{-}56& 31&\mbox{-}29&  0& 12&\mbox{-}26 \\
\mbox{-}5& \mbox{-}1& 22& 86&  3& 71& 39&  1&\mbox{-}22& \mbox{-}1& 31
&\mbox{-}18&  4
&  0&\mbox{-}17& 10& \mbox{-}9&  3& \mbox{-}1& \mbox{-}7 \\
\mbox{-}3&  0& \mbox{-}5&\mbox{-}37&  0&\mbox{-}29&\mbox{-}15& 11&  9&  1&\mbox{-}13
&  4&  0&  0&  8& \mbox{-}4&  4& \mbox{-}6&  1&  1 \\
 1&  0&  3& 26&  0& 22& 11&  4&\mbox{-}15&  0& 11&  0&  0&  0& \mbox{-}4&  0&  0&  1
&  0& \mbox{-}3 \\
12&  4&\mbox{-}15& \mbox{-}9&  0& \mbox{-}6& \mbox{-}1& \mbox{-}2& \mbox{-}1& \mbox{-}8
& \mbox{-}1& 16& \mbox{-}4&  0&  3& \mbox{-}8&  4&  3& \mbox{-}4& \mbox{-}1 \\
\mbox{-}3&  0& 11& 45&  0& 37& 19&  7&\mbox{-}15&  1& 19& \mbox{-}4&  0&  0& \mbox{-}8
&  4& \mbox{-}4&  0& \mbox{-}1& \mbox{-}7 \\
48& 16&\mbox{-}64&\mbox{-}63&  0&\mbox{-}47&\mbox{-}24& \mbox{-}7&  5& \mbox{-}7
&\mbox{-}18& 60&\mbox{-}16&  0& 16&\mbox{-}28& 15&  3& \mbox{-}3&  4 \\
16&  9& 10&145&  0&116& 55& 11&\mbox{-}38&  6& 49&  3& \mbox{-}9& \mbox{-}4&\mbox{-}26
& 10&\mbox{-}11&  0&  6&\mbox{-}11 \\
21&  8&\mbox{-}39&\mbox{-}93&  3&\mbox{-}75&\mbox{-}36& \mbox{-}9& 21& \mbox{-}3
&\mbox{-}33& 24& \mbox{-}3&\mbox{-}12& 18&\mbox{-}15& 12&  0& \mbox{-}3&  9 \\
\mbox{-}3&  0& \mbox{-}3&\mbox{-}18&  0&\mbox{-}12& \mbox{-}6&  0&  3&  0& \mbox{-}6
&  3&  0&  3&  6& \mbox{-}3&  3&  0&  0&  0 \\
\mbox{-}3&  1& 16& 57& \mbox{-}3& 41& 17&  5& \mbox{-}7&  3& 18&\mbox{-}11& \mbox{-}4
&  0&\mbox{-}11& 19&\mbox{-}11& \mbox{-}3&  3& \mbox{-}2 \\
 4&  4&\mbox{-}10&\mbox{-}18&  0&\mbox{-}12& \mbox{-}6& \mbox{-}3&  0&  0& \mbox{-}6
& 10& \mbox{-}4& \mbox{-}4&  6& \mbox{-}3&  6&  3&  0&  0 \\
15&  4&\mbox{-}24&\mbox{-}27&  0&\mbox{-}18& \mbox{-}7&\mbox{-}11& \mbox{-}4& \mbox{-}8
& \mbox{-}7& 25& \mbox{-}4&  0&  9&\mbox{-}17& 13& 12& \mbox{-}4& \mbox{-}1 \\
 1&  0&  3& 26&  0& 22& 11&  1&\mbox{-}15&  0& 11&  0&  0&  0& \mbox{-}4&  0&  0&  4
&  0& \mbox{-}3 \\
18&  4&\mbox{-}36&\mbox{-}95&  0&\mbox{-}77&\mbox{-}42& \mbox{-}7& 27&  3&\mbox{-}38
& 23& \mbox{-}4&  0& 18&\mbox{-}15& 11& \mbox{-}3&  5& 13 \\
22&  9& 10&177&  0&142& 68& 12&\mbox{-}53&  6& 62&  3& \mbox{-}9&  0&\mbox{-}32& 11
&\mbox{-}13&  1&  6&\mbox{-}11
\end{array} \right]$}
\]
\end{example}

\begin{wrapfigure}[8]{r}{3.3in}
\vspace{-4mm}
\epsfig{figure=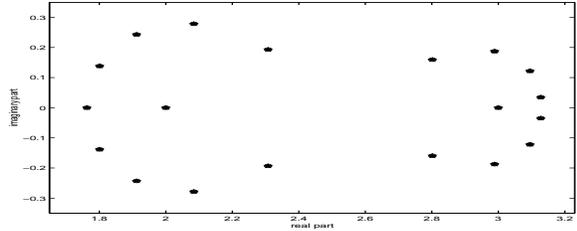,height=1.2in,width=3.0in}
\vspace{-4mm}
\caption{\footnotesize Eigenvalue clusters produced by Matlab} 
\label{fig:2clusters}
\end{wrapfigure}
has two exact eigenvalues ~$\la_1 \,=\, 2.0$~ and ~$\la_2\,=\, 3.0$~ with Segre
characteristics ~$\{9,1\}$~ and ~$\{8,2\}$~ respectively.
~Under round-off perturbation in the magnitude of machine precision
($\approx 2.2 \times 10^{\mns 16}$), Matlab outputs eigenvalues in two
noticeable clusters show in Figure \ref{fig:2clusters}.
~The arithmetic means of the two clusters are as follows

\vspace{-3mm} \scriptsize
\begin{verbatim}
                        means                exact eigenvalues    cluster condition number
        left cluster:   1.99724665369002     2.000000000000000         6.50e+012
        right cluster:  3.00275334630999     3.000000000000000         6.48e+012
\end{verbatim} \normalsize \vspace{-3mm}

From these results, we can see that only 3 correct digits are obtained by grouping.
~In contrast, our iterative method, which will be presented in \S \ref{s:ctrip},
converges on the two eigentriplets accurately
and attains 14 correct digits on the two eigenvalues.

%\vspace{-5mm}
\scriptsize
\begin{verbatim}

          Computed eigenvalues             2.00000000000004     3.00000000000003 
          ---------------------------------------------------------------------------
          forward error                    4.00e-15             3.02e-14
          backward error                   1.65e-17             5.77e-17
          staircase condition number       3.45e+07             5.33e+05
\end{verbatim}
\normalsize
\vspace{-3mm}

The cluster condition numbers are over
~$6 \times 10^{12}$~ and
the staircase condition numbers are substantially smaller
(~$3\times 10^7$ ~and ~$5 \times 10^5$).
~From the examples we have tested, computing 
staircase eigentriplet appears to be always more accurate than 
grouping clusters.

\vspace{-5mm}
\section{Computing a staircase eigentriplet with a known
structure} \label{s:ctrip}
\vspace{-5mm}

In this section we present the method for computing a numerical
unitary-staircase eigentriplet under the assumption that the Weyr
characteristic ~$\{m_1 \geq m_2 \ge \cdots\}$~
is known for an ~$m$-fold eigenvalue ~$\la$~ that is approximated by
~$\hat{\la}$.
~An algorithm for computing
the required Weyr characteristic and initial approximations to the
eigenvalues will be given in the next section (\S \ref{s:compstru}).
~There are two steps in calculating the staircase eigentriplet ~$(\la,U,S)$:
~First find an initial staircase eigentriplet
~$(\hat{\la}, \hat{U}, \hat{S})$,
~then the Gauss-Newton iteration is applied to refine the eigentriplet
until a desired accuracy is attained.

The QR decomposition and its updating/downdating will be used extensively.
~When a row is deleted from a matrix ~$B$~ to form a new matrix ~$\check{B}$,
~finding a QR decomposition of ~$\check{B}$~ from an existing
QR decomposition of ~$B$~ is called a QR downdating.
~Conversely, ~computing the QR decomposition after inserting a row
called a QR updating.
~QR updating and downdating are standard techniques in matrix computation
\cite[\S 12.5.3]{gvl}
requiring ~$O(m^2)$~ flops.

\vspace{-5mm}
\subsection{Computing the initial staircase eigentriplet}
\vspace{-5mm}

When ~$\hla \approx \la$~ is available with known multiplicity
~$m$~ and nonzero Weyr characteristic ~$m_1 \geq \cdots \geq m_k$, ~we
need an initial approximation ~$(\hla,\hU,\hS)$~ to the solution of
equations (\ref{auus}). 
~Write 
~$U=[U_1,\cdots,U_k]$~ with ~$U_j =
\blb \bdu_{\mu_{_{j\mns 1}}\pls 1}, \cdots, \bdu_{\mu_{_j}} \brb$. 
~From the uniqueness in Theorem \ref{unisc},  
each column ~$\bdu_{\mu_i \pls j}$ ~of ~$U_{i\pls 1}$ ~along with the ~$j$-th
column of ~$S$ ~is the unique
solution to the homogeneous system
\begin{equation} \label{inisys}
\left\{ \begin{array}{rcl}
(A-\la I) \, \bdu_{\mu_i\pls j} - U_1 S_{1,i\pls 1} \bde_j - \cdots -
U_i S_{i,i\pls 1} \bde_j & = & \bdo \\
\blb \bdu_1,\cdots, \bdu_{\mu_i\pls j\mns 1} \brb^\h \bdu_{\mu_i\pls j} & = & 
\bdo \\
\blb \bdb_{\mu_i\pls j \pls 1}, \cdots, \bdb_{\mu_{i\pls 1}}\brb^\h 
\bdu_{\mu_i\pls j} & =&   \bdo
\end{array} \right.
\end{equation}
up to a unit multiple 
for ~$i=0,\cdots,k\lmns 1$ ~and ~$j=1,\cdots,m_{i\pls 1}$ 
~where ~$S$~ is as in (\ref{sl}).
~Consequently, the vector ~$\bdz_{\mu_i\pls j}$ ~consists of components
~$\bdu_{\mu_i\pls j}$, $S_{1,i\pls 1} \bde_j, \cdots, S_{i,i\pls 1} \bde_j$
~spans the one-dimensional kernel of the matrix
\begin{equation} \label{Blj}
G_{i,j} =
\left[ \mbox{\scriptsize $\begin{array}{clc}
\blb \bdb_{\mu_{i}\pls j\pls 1},\cdots,
\bdb_{\mu_{i\pls 1}}\brb^\h &&  \\
A - \la I && -U_1,\cdots,-U_{i} \\
\blb \bdu_{1},\cdots,\bdu_{\mu_i\pls j\mns 1} \brb^\h &&
\end{array}$} \right],
~~~ i = 0, \cdots, k-1,  ~~ j=1, \cdots, m_{i\mns 1}
\end{equation}
Let ~$Q_{ij} \,R_{ij}$ ~be the QR decomposition of ~$G_{i,j}$.
~Then the vector ~$\bdz_{\mu_i\pls j}$ ~can be computed by a simple
inverse iteration \cite{li-zeng-rr} on ~$R = R_{ij}$
\begin{equation} \label{nvec}
\left\{ \begin{array}{l}
\mbox{set ~$\bdz_0$ ~as a random vector} \\
\mbox{for ~$j=1,2,\cdots$ ~do} \\
~~\left\lfloor \begin{array}{l}
\mbox{solve} ~R^\h \bdx = \bdz_{j\mns 1} \\
\mbox{solve} ~R \bdy = \bdx ~~\mbox{and set} ~\bdz_j = \bdy/\|\bdy\|_2
\end{array} \right.
\end{array} \right.
\end{equation}
After ~$\bdz_{mu_i\pls j}$ ~is computed from ~$G_{ij} = Q_{ij}R_{ij}$, 
~the next vector ~$\bdz_{mu_i\pls j\pls 1}$ ~will be computed from
~$G_{i,j\pls 1}$ ~which comes from deleting ~$\bdb_{mu_i\pls j\pls 1}$
~from the top row of ~$G_{ij}$ ~and inserting ~$\bdu_{\mu_i\pls j}$
~at the bottom.
~Namely, the QR decomposition of ~$G_{i,j\pls 1}$ ~is obtained from that of
~$G_{ij}$ ~via a QR updating and a QR downdating. 

In summary, computing the initial staircase eigentriplet
~$(\hla, \hU, \hS)$~ is
a process consisting of repeated QR updating/downdating and
consecutive applications of inverse iteration (\ref{nvec}), as outlined
in the following pseudo-code.

\vspace{-5mm}
\begin{itemize} \parskip-0.5mm
\item[] {\bf Algorithm} {\sc InitialEigentriplet}%
\index{Algorithm!{\sc InitialEigentriplet}}
\item[] Input: \ matrix ~$A$, \
Weyr char. ~$\{ m_1 \geq \cdots \ge m_k\}$, \
initial eigenvalue ~$\la = \hla$
\begin{itemize}
\item get random vectors ~$\bdb_1,\cdots,\bdb_m$~ and
QR decomposition of ~$A - \la I$ 
\item for ~$i=0,1,\cdots,k-1$ do
\item[] $\left\lfloor \begin{array}{l}
\mbox{Update the QR decomposition ~$G_{i1} = Q_{i1} R_{i1}$}\\
\mbox{for ~$j=1,2,\cdots,m_{i\pls 1}$~ do} \\
~\left\lfloor \begin{array}{l}
\mbox{apply iteration (\ref{nvec}) on ~$R_{ij}$ ~to find a
numerical null vector ~$\bdz$}\\
\mbox{extract ~$\hat{\bdu}_{\mu_{i}\pls j}$, $\hS_{1,i\pls 1}\bde_j, \cdots,
\hS_{i,i\pls 1}\bde_j$ ~from ~$\bdz$}\\
\mbox{get ~$G_{i,j\pls 1} = Q_{i,j\pls 1} R_{i,j\pls 1}$
~by QR downdating/updating on ~$G_{ij} = Q_{ij}R_{ij}$}
\end{array} \right. 
\end{array} \right.
$
%\begin{itemize}
%\item Update the QR decomposition ~$G_{i1} =
%Q_{i1} R_{i1}$~
%\item for ~$j=1,2,\cdots,m_{i\pls 1}$~ do
%\begin{itemize}
%\item apply iteration (\ref{nvec}) on ~$R_{ij}$ ~to find an
%approximate null vector ~$\bdz$
%\item extract ~$\hat{\bdu}_{\mu_{i}\pls j}$, $\hS_{1,i\pls 1}\bde_j, \cdots,
%\hS_{i,i\pls 1}\bde_j$ ~from ~$\bdz$
%\item 
%obtain ~$G_{i,j\pls 1} = Q_{i,j\pls 1} R_{i,j\pls 1}$
%~by QR downdating/updating on ~$G_{ij} = Q_{ij}R_{ij}$~
%\end{itemize}
%\item[] end do
%\end{itemize}
%\item[] end do
\end{itemize}
\item[] Output ~$\hU$, ~$\hS$
\end{itemize}
\vspace{-3mm}

{\bf Remark:} ~Computing a staircase form from a given eigenvalue
was proposed by Kublanovskaya \cite{kubl} in 1968.
~Ruhe \cite{ruhe-70-bit} improved the Kublanovskaya Algorithm
in 1970 by employing singular value decomposition (SVD) for determining
the numerical rank and kernel.
~Due to successive SVD computation, the original Kublanovskaya-Ruhe approach
leads to an ~$O(n^4)$ ~algorithm \cite{beelen-vandooren} in the worst case
senerio. 
~Further improvement has been proposed in 
\cite{beelen-vandooren,golub-wilkinson} that reduce 
the complexity to ~$O(n^3)$. 
~Our Algorithm {\sc InitialEigentriplet} can be considered a new improvement
from Kublanovskaya-Ruhe Algorithm. 
~The novelty of our algorithm includes (a) the nullity-one homogeneous system
(\ref{inisys}); ~(b) employing an efficient null-vector finder (\ref{nvec}) to 
replace the costly SVD; and ~(c) successive QR updating/downdating. 
~As a result, Algorithm {\sc InitialEigentriplet} is of complexity ~$O(n^3)$
~and fits our specific need in satisfying the second constraint in 
(\ref{auus}).
~Furthermore, our computation of staircase form goes further with a
refinement step using the Gauss-Newton iteration in the following section.
\hfill {\LARGE $\Box$}

\vspace{-3mm}
\subsection{Iterative refinement for a staircase eigentriplet}
\label{s:itref}
\vspace{-5mm}

The initial eigentriplet ~$(\hla,\hU,\hS)$~ produced by
Algorithm {\sc InitialEigentriplet} (or by existing variations of the
Kublanovskaya Algorithm) may not be accurate enough.
~One of the main features of our algorithm is an iterative refinement strategy 
for ensuring the highest achievable accuracy in computing the staircase 
eigentriplet. 
~We elaborate the process in the following.

Since ~$(\hla,\hU,\hS)$~ approximately satisfies (\ref{inisys}),
this eigentriplet is an approximate solution to (\ref{ayys}) for
~$\blb \bdc_1, \cdots, \bdc_m \brb =
\blb \hat{\bdu}_1, \cdots, \hat{\bdu}_m \brb  =\hU$.
~Using these ~$\bdc_i$'s in (\ref{ayys}) and (\ref{cF}), we apply
the Gauss-Newton iteration for ~$i=0,1,\cdots$,
\begin{equation} \label{cgnrfn}
\big\llbracket \la^{(i\pls 1)}, Y^{(i\pls 1)}, S^{(i\pls 1)}\big\rrbracket ~=~
\big\llbracket \la^{(i)}, Y^{(i)}, S^{(i)} \big\rrbracket -
J\big( \la^{(i)}, Y^{(i)}, S^{(i)} \big)^+
\bdf\big( \la^{(i)}, Y^{(i)}, S^{(i)} \big)
\end{equation}
with initial iterate
~$\big\llbracket \la^{(0)}, Y^{(0)}, S^{(0)} \big\rrbracket ~=~
\llbracket \hla,\hU,\hS\rrbracket$.
~Here again, 
~$\llbracket \la^{(i)}, Y^{(i)}, S^{(i)} \rrbracket$~ denotes the vector
form of ~$\left( \la^{(i)}, Y^{(i)}, S^{(i)} \right)$~ according to the rule
given in (\ref{triporder}).

Let ~$(\la,Y,S)$ ~be a least squares solution of ~$(\ref{ayys})$ ~with 
sufficiently small residual, or equivalently ~$A$ ~is close to a matrix
~$\hat{A}$ ~having ~$(\la,Y,S)$ ~as its exact eigentriplet.
~Then the Jacobian ~$J(\cdot,\cdot,\cdot)$ ~is injective
by Theorem~\ref{t:cplxreg}, ensuring the Gauss-Newton iteration (\ref{cgnrfn}) 
to converges to ~$(\la,Y,S)$ ~locally. 
~This ~$(\la,Y,S)$ ~is a numerical staircase eigentriplet, a 
{\em unitary}-staircase eigentriplet can be obtained by an orthogonalization
and an extra step of refinement.
~Specifically, Let ~$Y = U\,R$ ~be 
the ``economic'' QR decomposition with 
~$U=\blb \bdu_1,\cdots,\bdu_m \brb$. 
~Partitioning ~$U = \blb U_1,\cdots,U_k\brb$~ the same way
as ~$Y = \blb Y_1,\cdots,Y_k\brb$, ~it is straightforward to verify
that ~$\cR{\blb U_1,\cdots,U_l\brb} = \cN{(A-\la I)^l}$ ~for
~$l = 1, \cdots, k$ ~if ~$A$ ~has ~$\la$ as an exact eigenvalue with
Weyr characteristic ~$\{m_j\}$, ~and ~$U^\h A U - \la I$ ~is the corresponding
nilpotent staircase form. 
~Furthermore, by resetting ~$S = U^\h A U - \la I$, 
~$\bdc_1 = \bdb_1 = \bdu_1$, ~$\cdots$, ~$\bdc_m= \bdb_m = \bdu_m$, 
the equations in ~(\ref{ayys}) are satisfied including the auxilliary 
equations.

In actual computation with the empirical data matrix ~$A$, ~a small error may 
emerge during the reorthogonalization process. 
~This error can easily be eliminated by one extra step of
refinement via the Gauss-Newton iteration starting from the new eigentriplet
~$(\la,U,S)$. 

\vspace{-5mm}
\begin{itemize} \parskip-0.5mm
\item[] {\bf Algorithm} {\sc EigentripletRefine}%
\index{Algorithm!{\sc EigentripletRefine}}
\item[] Input: \ Initial approximate unitary-staircase eigentriplet
~$(\hla,\hU,\hS)$, ~tolerance ~$\dl > 0$
\begin{itemize}
\item Set ~$\left( \la^{(0)}, Y^{(0)}, S^{(0)} \right) ~=~ (\hla,\hU,\hS)$,
~and ~$\blb \bdc_1, \cdots, \bdc_m \brb ~=~ \hU$
\item For ~$i = 1, 2, \cdots$~ do
\item[] ~$\left\lfloor \begin{array}{l}
\mbox{Solve
~\( J\left( \la^{(i\mns 1)}, Y^{(i\mns 1)}, S^{(i\mns 1)}\right) \bdz ~=~
\bdf\left( \la^{(i\mns 1)}, Y^{(i\mns 1)}, S^{(i\mns 1)}\right) \)~
for ~$\bdz$} \\
\mbox{Set ~$\left[\la^{(i)}, Y^{(i)}, S^{(i)}\right]  ~=~
\left[ \la^{(i\mns 1)}, Y^{(i\mns 1)}, S^{(i\mns 1)} \right] - \bdz$}\\
\mbox{If ~$\|\bdz\|_2 < \dl $~, ~then set
~$(\la,Y,S) ~=~ \left(\la^{(i)}, Y^{(i)}, S^{(i)}\right)$~
and break the loop} \end{array} \right.$
%\begin{itemize}
%\item Solve
%~\( J\left( \la^{(i\mns 1)}, Y^{(i\mns 1)}, S^{(i\mns 1)}\right) \bdz ~=~
%\bdf\left( \la^{(i\mns 1)}, Y^{(i\mns 1)}, S^{(i\mns 1)}\right) \)~
%for ~$\bdz$
%\item Set ~$\left[\la^{(i)}, Y^{(i)}, S^{(i)}\right]  ~=~
%\left[ \la^{(i\mns 1)}, Y^{(i\mns 1)}, S^{(i\mns 1)} \right] - \bdz$
%\item If ~$\|\bdz\|_2 < \dl $~, ~then set
%~$(\la,Y,S) ~=~ \left(\la^{(i)}, Y^{(i)}, S^{(i)}\right)$~
%and break the loop
%\end{itemize}
%\item[] end do
\item Economic QR decomposition ~$Y ~=~ U R$~ and set ~$S ~=~ U^\h A U - \la I$
\item If ~$R \approx I$, ~exit. ~Otherwise, 
set ~$(\hla,\hU,\hS) ~=~ (\la,U,S)$, 
~reset ~$\bdb_1,\cdots,\bdb_m$ ~and ~$\bdc_1,\cdots,\bdc_m$, 
and repeat the algorithm
\end{itemize}
\item[] Output ~$(\la,U,S)$
\end{itemize}
\vspace{-5mm}

To carry out the iterative refinement (\ref{cgnrfn}), a QR
decomposition ~$ J\left( \la^{(i)}, Y^{(i)}, S^{(i)} \right) = Q_i R_i$
~is required at every iteration step.
~A straightforward QR decomposition costs
~$O\left( (mn)^3 \right)$~ flops, which can be substantially reduced by
taking the structure of the Jacobian into account.
~Using the Weyr characteristic ~$\{3,2,1\}$ ~as an example, the nilpotent
staircase matrix ~$S$ ~is shown in (\ref{S321}). 
~the Jacobian ~$J(\la,Y,S)$ ~of ~$\bdf(\la,Y,S)$~ 
with a re-arrangement of columns and rows ~$P_1$~ and ~$P_2$
~becomes
\[
 \begin{array}{l} \mbox{\normalsize $P_1J(\la,Y,S)P_2 ~~=~~$} \\ \\
\left[ \mbox{\tiny \begin{tabular}{llllllrccc}\cline{1-1}
\multicolumn{1}{|l|}{$C_6^\h$
%\begin{array}{c} \bdc_1^\h \\ \vdots \\ \bdc_6^\h \end{array}$
}
& & & & & & & & & \\
\multicolumn{1}{|c|}{$A-\la I$} &$\mns s_{56}I$ &$\mns s_{46}I$ &
$\mns s_{36}I$ & $\mns s_{26}I$ &$\mns s_{16}I$ &
$\mns \bdy_6$ &  & & \hspace{-4mm} $\mns [\bdy_1, \cdots,\bdy_5]$ \\ \cline{1-2}
\multicolumn{1}{c}{ } &
\multicolumn{1}{|l|}{$C_5^\h$
%$\begin{array}{c} \bdc_1^\h \\ \vdots \\ \bdc_5^\h \end{array}$
} & & & &  & & & & \\
\multicolumn{1}{c|}{ } & \multicolumn{1}{|l|}{$A-\la I$}
& & $\mns s_{35}I$ & $\mns s_{25}I$ &
$\mns s_{15} I$   & $\mns \bdy_5$ & &
\hspace{-4mm} $\mns [\bdy_1,\bdy_2, \bdy_3]$ &\\ \cline{2-3}
& \multicolumn{1}{l}{ } & \multicolumn{1}{|l|}{$\bdb_{5}^\h$} & & & & & & &   \\
\multicolumn{2}{c|}{ }  & \multicolumn{1}{|l|}{$C_4^\h$
%$\begin{array}{l} \bdc_1^\h \\ \vdots \\ \bdc_4^\h \end{array}$
} & & & & & & &  \\
\multicolumn{2}{c|}{ }  &  \multicolumn{1}{|l|}{$A-\la I$}
& $\mns s_{34}I$ & $\mns s_{24}I$ &
$\mns s_{14} I$   & $\mns \bdy_4$ & $\mns [\bdy_1,\bdy_2, \bdy_3]$ & & \\
\cline{3-4}
\multicolumn{3}{c|}{ }  & \multicolumn{1}{|l|}{$C_3^\h$
%$\begin{array}{l} \bdc_1^\h \\ \bdc_2^\h \\ \bdc_3^\h \end{array}$
}  & & & & & & \\
\multicolumn{3}{c|}{ } & \multicolumn{1}{|l|}{
$A-\la I$} & & & $\mns \bdy_3$ & & &  \\ \cline{4-5}
\multicolumn{4}{c|}{ } & \multicolumn{1}{|l|}{$\bdb_{3}^\h$} & & & & &  \\
\multicolumn{4}{c|}{ } & \multicolumn{1}{|l|}{$C_2^\h$
%$\begin{array}{l} \bdc_1^\h \\ \bdc_2^\h \end{array}$
} & & & & &  \\
\multicolumn{4}{c|}{ } & \multicolumn{1}{|l|}{
$A-\la I$} & & $\mns \bdy_2$ & & &  \\ \cline{5-6}
\multicolumn{5}{c|}{ } & \multicolumn{1}{|l|}{$\bdb_{2}^\h$} & & & &  \\
\multicolumn{5}{c|}{ } & \multicolumn{1}{|l|}{$\bdb_{3}^\h$} & & & &  \\
\multicolumn{5}{c|}{ } & \multicolumn{1}{|l|}{$C_1^\h$} & & & & \\
\multicolumn{5}{c|}{ } & \multicolumn{1}{|l|}{$A-\la I$}
& $\mns \bdy_1$ & & &  \\ \cline{6-6}
\end{tabular}} \right] \end{array}
\]
where the blocks ~$C_j = \blb \bdc_1,\cdots,\bdc_j\brb$ ~for ~$j=1,2,\cdots,
6$.
~Without loss of generality, we can assume ~$A$~ is in Hessenberg form.
~Then ~$J(\la,Y,S)$~ is near triangular and requires only ~$O(m^3n^2)$~
flops for its QR decomposition.
~The backward substitution requires ~$O(m^2 n^2)$ ~flops.

\vspace{-5mm}
\subsection{Converting a staircase form to Jordan decomposition}
\label{sec:sf2jd}
\vspace{-5mm}

After finding a staircase decomposition ~$A = UTU^\h$,
~the Jordan Canonical Form of ~$A$~ is available by conjugating the
Weyr characteristic.
~In the cases where the Jordan decomposition is in demand,
a method for converting the staircase decomposition to the Jordan
decomposition is proposed by Kublanovskaya \cite{kubl} which
primarily involves {\em non-unitary}
similarity transformations.
~Detailed procedures
can also be found in \cite{golub-wilkinson,kagstrom-ruhe}.

Alternatively, we may calculate the Jordan decomposition
~$A X = X J$~ by converting unitary-staircase eigentriplets
~$(\la_i,U_i,S_i)$
~to Jordan decompositions ~$A (U_i G_i) \,=\, (U_i G_i) J_i$~ for
~$i=1,\cdots,k$~ using Kublanovskaya's algorithm.

\vspace{-5mm}
\subsection{Numerical examples for computing the staircase form}
\label{sec:numres1}
\vspace{-5mm}

Algorithm {\sc InitialEigentriplet} combined with Algorithm
{\sc EigentripletRefine} forms a stand-alone algorithm for
computing a staircase/Jordan decomposition,
assuming an initial approximation to a multiple eigenvalue together with its
Segre/Weyr characteristics are available by other means.
~This combination is implemented as a Matlab module {\sc EigTrip}.
We shall present a method for computing the Jordan structure in
\S \ref{s:compstru}.

A previous algorithm for computing a staircase form with 
known Segre characteristic via a minimization process is constructed by 
Lippert and Edelman \cite{lip-edel-00} and implemented as a Matlab module
{\sc sgmin}.
~The iteration implemented in {\sc sgmin} converges in many cases,
including difficult test matrices such as the Frank matrix.
~As we shall show below, our algorithm provides a substantial improvement
over {\sc sgmin} particularly on cases where the cluster condition numbers 
(\ref{clcond})
are large but the staircase condition numbers stay moderate.
~We list the comparisons on accuracy only.

\setcounter{example}{0}

\vspace{-3mm}
\begin{example} \hspace{-2.5mm}' \em
~~We test both {\sc sgmin} and our {\sc EigTrip} on the matrix ~$A$~
given in Example \ref{e:sens} in \S \ref{sec:sens}, starting from eigenvalue
approximation
~$\hla_1 = 1.999$~ and ~$\hla_2= 2.999$~ with
given Segre characteristics ~$\{9,1\}$~ and ~$\{8,2\}$~ respectively.
~The code {\sc sgmin} improves the eigenvalue accuracy by one and four digits 
respectively.
~In contrast, our {\sc EigTrip} obtains an accuracy near the
machine precision on both eigenvalues, as shown in Table \ref{tbl:ex1}.
\end{example}

\begin{table}[ht]  \small
\begin{center}
\begin{tabular}{|l|l|c|l|c|} \hline
& \multicolumn{2}{c|}{from ~$\hla_1=1.999$} & \multicolumn{2}{c|}{from ~$\hla_2=2.999$}
\\ \cline{2-5}
& computed & backward & computed & backward \\
& eigenvalue & error & eigenvalue & error \\ \hline
cluster mean & {\bf 1.99}724665369002  &  --- & {\bf 3.00}275334630999 & --- \\
{\sc sgmin} &  {\bf 1.9999}1878946447  & 1.004e-008  & {\bf 2.9999999}1118127
& 6.895e-010 \\
{\sc EigTrip} & {\bf 1.9999999999999}8 & 3.270e-017  & {\bf 3.00000000000000}3
& 4.673e-017 \\ \hline
\end{tabular}
\end{center} \vspace{-3mm}
\caption{\footnotesize Accuracy comparison for Example 1} \label{tbl:ex1}
\end{table}

\vspace{-3mm}
\begin{example} \label{e:50m}
\em ~~We construct a ~$50\times 50$~ matrix having known
multiple eigenvalues ~$\la = 1.0, ~2.0$~ and ~$3.0$~ with Segre characteristics
~$\{10,5,3,2\}$, ~$\{8,4,3\}$~ and ~$\{4,1\}$~ respectively, together with
ten simple eigenvalues randomly generated in the
box ~$[-3,3]\times [-3,3]$.
~Both {\sc sgmin} and {\sc EigTrip} start at initial approximations
~$\la_1^{(0)} = 0.99$, ~$\la_2^{(0)} = 1.99$, and ~$\la_3^{(0)} = 2.99$.
~The results of the iterations are listed in Table \ref{tbl:50m}, in which
forward errors are ~$|\la_j - \hat{\la}_j|$~ for each computed eigenvalue
~$\hat{\la}_j$, ~$j=1,2,3$, and the backward errors are the residual 
(\ref{resdef})
for each eigentriplet.
 \end{example}

\begin{table}[ht]
\begin{center}
\begin{tabular}{|l|cc|cc|cc|} \hline
& \multicolumn{2}{c|}{at ~$\la = 1.0$} & \multicolumn{2}{c|}{at ~$\la = 2.0$}
& \multicolumn{2}{c|}{at ~$\la = 3.0$}  \\ \cline{2-7}
& forward & backward & forward & backward & forward & backward  \\
& error & error & error & error & error &  error \\ \hline
{\sc sgmin} & 2.29e-008& 8.46e-007& 5.01e-008 & 9.42e-007 & 1.03e-009 & 3.15e-008\\ 
\hline
{\sc EigTrip} & 2.22e-016 & 1.16e-015& 0 & 1.89e-016 & 8.88e-016& 1.23e-016  \\ 
\hline
\end{tabular}
\end{center} \vspace{-3mm}
\caption{\footnotesize
Accuracy comparison for Example \ref{e:50m}} \label{tbl:50m}
\end{table}

The results show that our algorithm is capable of calculating eigenvalues to
the accuracy near machine precision (16 digits).
~For each approximate eigentriplet ~$(\la,Y,S)$~ of matrix ~$A$,
~the residual ~$\rho$~ is defined in (\ref{resdef}).
~By (\ref{resdis}), with relative distance up to ~$\rho$~ from ~$A$,
~there is a nearby matrix ~$\hat{A}$~
for which ~$(\la,Y,S)$~ is an exact eigentriplet.

\vspace{-3mm}
\begin{example} {\bf (Frank matrix)}  \em
\cite{chatelin-fraysse,golub-wilkinson,kagstrom-ruhe,lippert-edelman,%
ruhe-70-bit,tref-emb,wilkinson}:
~This is a classical test matrix given in a Hessenberg form
~$F = \left( f_{ij} \right)$,  ~with
~$f_{ij} = n+1-\max\{i,j\}$ ~for $j\ge i-1$
~and ~$f_{ij}=0$ ~otherwise.
~Frank matrix has no multiple eigenvalues.
~However, its small eigenvalues are ill-conditioned measured by the
standard eigenvalue condition number \cite{gvl},
as shown in the following table.  \end{example}

\scriptsize
\begin{center} \begin{tabular}{||r|r||r|r||r|r||} \hline\hline
\multicolumn{6}{||c||}{\small
Eigenvalues and condition numbers of ~$12\times 12$~ Frank matrix} \\ \hline
\small Eigenvalue & \small condition & \small Eigenvalue &
\small condition \ \ \ & \small Eigenvalue & \small condition \ \ \
\\ \hline\hline
 32.22889 &    8.5  &  3.51186 &         34.1 & 0.143647 &  611065747.8 \\
 20.19899 &   16.2  &  1.55399 &       1512.5 & 0.081228 & 2377632497.8 \\
 12.31108 &    9.0  & 0.64351 &    1371441.3 & 0.049507 & \ \ \ 3418376227.8 \\
  6.96153 &   24.1  & 0.28475 &   53007100.5 & 0.031028 & \ \ \ 1600156877.4
\\ \hline\hline
\end{tabular} \vspace{-3mm} \end{center}
\normalsize

Clearly, Frank matrix is near matrices which possess multiple eigenvalues
near zero with nontrivial Jordan structures.
~Using an initial eigenvalue estimation near zero and Segre characteristics
~$\{2\}$, ~$\{3\}$, ~$\{4\}$, ~$\{5\}$~ and ~$\{6\}$~ in consecutive tests,
our refinement algorithm {\sc EigTrip} produces five nearby matrices with an
eigenvalue of multiplicity 2, 3, 4, 5, and 6
respectively, as shown in the table below.

\small
\begin{center} \begin{tabular}{|l|cccll|} \hline
& \multicolumn{5}{|c|}{ \normalsize 5 nearby matrices with following features
respectively}  \\
\cline{2-6}
& \normalsize given & \normalsize \ \ computed \ \ &
\normalsize \ \ backward \ \
& \normalsize \ \ staircase \ \ & cluster \\
& \normalsize Segre ch. & \normalsize \ \ eigenvalue \ \ &
\normalsize \ \ error \ \
& \normalsize \ \ condition \ \ & condition \\ \hline
{\sc sgmin} & \{6\}   & 0.1870511240986754 & 6.34e-05 &      & 126.8\\
{\sc EigTrip} & \{6\} & 0.1870509025041315 & 6.34e-05 & 5.96 & \\ \hline
{\sc sgmin} & \{5\}   & 0.1076751260727581 & 1.90e-06 &      & 7689.2 \\
{\sc EigTrip} & \{5\} & 0.1076751114381528 & 1.90e-06 & 32.2 & \\ \hline
{\sc sgmin} & \{4\}   & 0.0701182985767899 & 6.12e-08 &      & 291589.8\\
{\sc EigTrip} & \{4\} & 0.0703019426541069 & 3.47e-08 & 447.4 & \\ \hline
{\sc sgmin} & \{3\}   & 0.0504328996330119 & 4.23e-l0 &      & 3666804.6 \\
{\sc EigTrip} & \{3\} & 0.0504338685708545 & 4.23e-10 & 11322.9 & \\ \hline
{\sc sgmin} & \{2\}   & 0.0305042120283680 & 9.87e-10 &      & 15192435.2\\
{\sc EigTrip} & \{2\} & 0.0386493437615946 & 3.45e-12 & 458607.1 & \\ \hline
\end{tabular} \end{center}
\normalsize

In other words, Frank matrix ~$F$~ resides within a relative distance
~$3.45 \times 10^{\mns 12}$~ from a matrix having a double eigenvalue,
or ~$4.23 \times 10^{\mns 10}$~ from a matrix having a triple eigenvalue,
etc.
~Notice that the cluster condition numbers (\ref{clcond}) in both cases
are quite high whereas the staircase condition numbers are small.
~It appears that our Algorithm {\sc EigTrip} substantially improves 
backward accuracy over {\sc sgmin}, 
particularly when cluster condition number is large.

\vspace{-5mm}
\section{Computing the numerical Jordan structure}
\label{s:compstru}
\vspace{-5mm}

In this section we present the theory and algorithm for computing
the structure of the numerical Jordan Canonical Form represented by
Segre and Weyr characteristics.

\vspace{-5mm}
\subsection{The minimal polynomial}
\vspace{-5mm}

As described in many textbooks on fundamental algebra (see, e.g.
\cite{artin}), given a linear operator ~$T\,:\; \cV \longrightarrow \cV$
~on a vector space ~$\cV$~ over a field ~$\cF$, ~one may view ~$\cV$~ as a
{\em module} over ~$\cF[t]$\index{$\cF[t]$: \ \ polynomial ring over $\cal F$}
~by a ``scalar'' product:
~$p(t)\bdv ~\equiv~ p(T) \bdv ~=~
a_n T^n(\bdv) + a_{n\mns 1} T^{n\mns 1}(\bdv) +
\cdots + a_1 T(\bdv) + a_0 \bdv$
~for ~$p(t) = a_n t^n + \cdots + a_1 t + a_0 \in \cF[t]$~ and ~$\bdv \in \cV$.
~For ~$\cF = \bdC$,
~$\cV = \bdC^n$, ~and ~$A\in \bdC^{n \tms n}$~ being the matrix representation
of ~$T$, ~we consider ~$\bdC^n$~ a module over the polynomial
ring ~$\bdC [t]$~ with scalar product ~$p(t)\bdv \equiv p(A)\bdv$~ for
~$p(t) \in \bdC[t]$~ and ~$\bdv \in \cV$.

A monic polynomial ~$p(t) \in \bdC[t]$~ is called an {\em annihilating
polynomial}\index{annihilating polynomial}\index{polynomial!annihilating}
for ~$\bdv \in \cV$ (with respect to $A$) ~if
~$p(t) \bdv \;(\,\equiv p(A)\bdv\,) = \bdo$.
~For a subspace ~$\cW \subseteq \bdC^n$, ~if ~$p(t)\bdv = \bdo$~ for all
~$\bdv \in \cW$,
~then ~$p(t)$~ is regarded as an annihilating polynomial for ~$\cW$.
~The polynomial with least degree among all the annihilating polynomials
for ~$\bdv$ (or subspace ~$\cW$)~ is called the
{\em minimal polynomial}\index{minimal
polynomial}\index{polynomial!minimal} for ~$\bdv$~ (or subspace ~$\cW$).
Note that every annihilating polynomial for ~$\bdv$~ (or subspace ~$\cW$) is
divisible by the minimal polynomial and obviously the minimal polynomial
for subspace ~$\cW$~ is divisible by any minimal polynomial for any vector
in ~$\cW$.
~If the minimal polynomial for a vector ~$\bdv \in \cW$~ coincides with the
minimal polynomial for ~$\cW$~ then ~$\bdv$~ is said to be a {\em regular
vector}\index{regular vector (of a subspace)} of ~$\cW$.

By the Fundamental Structure Theorem for modules over Euclidean
domain \cite{artin}, ~$\bdC^n$~ is a direct sum of cyclic submodules,
say ~$\bdC^n = \cW_1 \oplus \cdots \oplus \cW_k$,
~where for each ~$i = 1, \cdots, k$,
~$\cW_i$~ is a cyclic submodule
(a submodule spanned by one vector) invariant with respect to ~$A$~
and is isomorphic to ~$\bdC[t]/\left(p_i(t)\right)$~ with ~$p_i(t)$~
being the minimal polynomial for ~$\cW_i$.
~Moreover, each ~$p_i(t)$~ is divisible by ~$p_{i\pls 1}(t)$~ for
~$i=1,\cdots,k-1$, ~that is
~$p_k(t) \divs p_{k\mns 1}(t) \divs \cdots \divs p_1(t)$.
Here, for polynomial ~$h(t)$~ and ~$q(t)$, ~notation
~$h(t)\divs q(t)$~ stands for
~``$h(t)$~ divides ~$q(t)$''.

It follows from ~$p_k(t) \divs \cdots \divs p_1(t)$
~that
each ~$p_i(t)$ ~for ~$i=1,\cdots,k$~ can
be written in the form
\begin{equation} \label{piform}
p_i(t) ~~=~~ (t-\al_1)^{m_{i1}} \cdots (t-\al_l)^{m_{il}}
\end{equation}
for fixed ~$\al_1,\cdots,\al_l \in \bdC$, ~and ~$m_{1j} \geq m_{2j} \geq
\cdots \geq m_{kj} \ge 0$~ for ~$j=1,\cdots,l$.
\vspace{-4mm}
\begin{lemma} {\em \cite{artin}} \label{eJb}
~~For each ~$(t-\al_j)^{m_{ij}}$~ in {\em (\ref{piform})} with ~$m_{ij} > 0$~
where ~$i=1,\cdots,k$~ and ~$j=1,\cdots,l,$~ there is an
elementary Jordan block ~$J_{m_{ij}}(\al_j)$
~of order ~$m_{ij}$~ associated with eigenvalue ~$\al_j$
in the Jordan Canonical Form of ~$A$, ~and the Jordan Canonical Form
of ~$A$~ consists of all such elementary Jordan blocks.
\end{lemma}
\vspace{-4mm}

When a subspace ~$\cW \subset \bdC^n$~ is invariant with respect to
~$A$, ~the linear transformation ~$A$~ induces a linear map
~$\tilde{A}\,:\,\bdC^n/\cW \,\longrightarrow\,\bdC^n/\cW$
~given by ~$\tilde{A}(\bdv+\cW) = A\bdv + \cW$.
~All the concepts and statements on annihilating polynomials and
minimal polynomials introduced above for ~$\bdC^n$~ with
linear map ~$A\,:\,\bdC^n \rightarrow \bdC^n$~ can be repeated
for ~$\bdC^n/\cW$~ with linear map
~$\tilde{A}\,:\,\bdC^n/\cW \,\rightarrow\,\bdC^n/\cW$.
~For instance, ~$p(t) \in \bdC[t]$~ is the minimal polynomial for subspace
~$\tilde{\cU} \subset \bdC^n/\cW$~ if ~$p(t)$~ is the least degree
polynomial which annihilates all ~$\tilde{\bdu} \in \tilde{\cU}$, ~that is
~$p(t)\tilde{\bdu} = p(\tilde{A}) \tilde{\bdu} = \bdo$~ for all
~$\tilde{\bdu} \in \tilde{\cU}$.

\begin{wrapfigure}[6]{r}{3.5in}
\vspace{-4mm}
\centerline{
$
\begin{CD} \bdC^n/\cW @>\tilde{A}>> \bdC^n/\cW \\
@V{\sg}VV  @VV{\sg}V \\
\cW' @>B>> \cW'
\end{CD}
$}
\vspace{-0mm}
\caption{\footnotesize Commuting diagram} 
\label{fig:cd}
\end{wrapfigure}
When ~$\bdC^n/\cW$~ is isomorphic to a vector space ~$\cW'$~ over
~$\bdC$~ with isomorphism ~$\sg\,:\,\bdC^n/\cW \rightarrow \cW'$,
~then the linear map ~$\tilde{A}\,:\,\bdC^n/\cW \rightarrow \bdC^n/\cW$~
induces a linear map
~$B = \sg \circ \tilde{A} \circ \sg^{\mns 1} \,:\, \cW' \longrightarrow
\cW'$,
~making the diagram in Figure~\ref{fig:cd} commutes.
~That is, ~$B\circ \sg \;=\; \sg \circ \tilde{A}$.

\vspace{-4mm}
\begin{lemma} \ \label{mpln}
~~For any subspace ~$\tilde{\cU} \subset \bdC^n/\cW$, ~$p(t) \in \bdC[t]$~
is the minimal polynomial for ~$\tilde{\cU}$ ~with respect to ~$\tilde{A}$
~if and only if ~$p(t)$~ is the minimal polynomial of ~$\sg(\tilde{\cU})$~
with respect to ~$B$.
\end{lemma}
\vspace{-4mm}

\prf ~$B = \sg \circ \tilde{A} \circ \sg^{\mns 1}$ ~implies ~$B^m =
\sg \circ \tilde{A}^m \circ \sg^{\mns 1}$~ for any integer
~$m > 0$. 
~It follows that
~$g(B) = \sg \circ g(\tilde{A}) \circ \sg^{\mns 1}$
~for any ~$g(t) \in \bdC[t]$.
~Thus ~$g(B) \sg(\tilde{\bdu}) = \sg \circ g(\tilde{A}) \tilde{\bdu}$
~for ~$\tilde{\bdu} \in \tilde{\cU}$
~and 
\begin{equation} \label{lra}
 g(B) \sg(\tilde{\bdu}) ~=~ 0 ~~\Longleftrightarrow~~
g(\tilde{A})\, \tilde{\bdu} ~=~ 0.
\end{equation}
Let ~$p_1(t)$~ be the minimal polynomial for ~$\tilde{\cU}$~ (with
respect to ~$\tilde{A}$) and ~$p_2(t) $~ be the minimal polynomial
for ~$\sg(\tilde{\cU})$~ (with respect to ~$B$).
~Then, by (\ref{lra}), ~$p_1(\tilde{A}) \tilde{\bdu} = 0$~ implies
~$p_1(B) \sg(\tilde{\bdu}) = \bdo$ for all $\tilde{\bdu} \in \tilde{\cU}$.
~So, ~$p_1(t)$~ annihilates ~$\sg(\tilde{\cU})$~ and hence ~$p_2(t) \divs
p_1(t)$.
~By the same argument ~$p_1(t) \divs p_2(t)$, ~and the
assertion follows. \hfill {\LARGE $\Box$}

\vspace{-5mm}
\subsection{The Jordan structure via minimal polynomials}
\label{sec:mp}
\vspace{-5mm}

By Lemma \ref{eJb}, the first task in finding the Jordan structure 
of ~$A\,:\,\bdC^n \rightarrow \bdC^n$~ is to
identify the minimal polynomial ~$p_i(t)$~ for the corresponding
cyclic submodules ~$\cW_i$, ~$i=1,\cdots,k$~ in 
~$\bdC^n = \cW_1 \oplus \cdots \oplus \cW_k$,
followed by factorizing ~$p_i(t)$~ in the form given in (\ref{piform}).
~We must emphasize here that accurate factorization of ~$p_i(t)$ ~in
numerical computation used to be regarded as a difficult problem.
~However, the appearance of a newly developed numerical algorithm
{\sc MultRoot} \cite{zeng_multroot,zeng-05} for calculating multiple roots
and their multiplicities makes this problem well-posed and solvable.
~Consequently the structure of the Jordan Canonical Form can be determined 
accurately.

We shall begin by finding the minimal polynomial ~$p_1(t)$~ for
~$\cW_1$.
~From ~$\bdC^n = \cW_1 \oplus \cdots \oplus \cW_k$, every ~$\bdv\in\bdC^n$~
can be written in the form
~$\bdv ~=~ \bdv_1 + \cdots + \bdv_k$ ~where
~$\bdv_i \in \cW_i$ ~for ~$i = 1, \cdots, k$.
~Thus, by ~$p_k(t) \divs \cdots \divs p_1(t)$, 
~we have ~$p_1(t) \bdv = p_1(A) \bdv = p_1(A) \bdv_1 + \cdots + p_1(A)
\bdv_k = \bdo$,
~making ~$p_1(t)$~ the minimal polynomial for ~$\bdC^n$.
~Meanwhile, ~$p_1(t)$~ is the minimal polynomial for all
~$\bdv \in \bdC^n$~ except those ~$\bdv$'s ~for which ~$\bdv_1 = \bdo$.
~The exceptional set is of measure zero.
~Therefore almost every
~$\bdv \in \bdC^n$~ is a regular vector.
~In other words, vector ~$\bdv$ ~is regular with probability one 
if it is chosen at random as in \S \ref{s:minpoly}.

To find minimal polynomial ~$p_1(t)$, ~we choose a generic
vector ~$\bdx \in \bdC^n$~ and check the
dimensions of the Krylov subspaces
~$\spn{\bdx,A\bdx}$, ~$\spn{\bdx,A\bdx,A^2\bdx}$, 
~$\spn{\bdx,A\bdx,A^2\bdx,A^3\bdx}$, ~$\cdots$
consecutively to look for the first integer ~$j$~ where
~$\spn{\bdx,A\bdx,\cdots,A^j\bdx}$~ is of dimension ~$j$.
~For this ~$j$, ~let
~$c_0'\bdx + c_1'A\bdx + \cdots + c_j'A^j\bdx = \bdo$.
%\[ \blb \bdx, A\bdx, \cdots, A^j\bdx \brb
%\left( \mbox{\scriptsize
%$\begin{array}{c} c_0'\\ \vdots \\ c_j' \end{array}$} \right)
%= \bdo.
%\]
~Obviously, ~$c_j' \neq 0$~ and
\[ p_1(t) ~~=~~ t^j + c_{j\mns 1} t^{j\mns 1} + \cdots + c_0, \mbox{~~with~~}
c_i = c_i'/c_j, \;\;\; i = 1, \cdots, j\lmns 1
\]
can serve as the minimal polynomial of ~$\cW_1$.
~We then proceed to find the minimal polynomial ~$p_2(t)$~ for ~$\cW_2$.
~By the same argument given above along with the property
~$\bdC^n/\cW_1 \simeq \cW_2 \oplus \cdots \oplus \cW_k = \cW'$,
~$p_2(t)$~ is the minimal
polynomial for ~$\cW'$~ (by ~$p_k(t) \divs \cdots \divs p_1(t)$ 
~as well as the minimal polynomial for almost all ~$\bdv \in \cW'$.
~By Lemma \ref{mpln}, ~$p_2(t)$~ is the
minimal polynomial for ~$\bdC^n/\cW_1$~ (with respect to the induced
linear map ~$\tilde{A} \,:\, \bdC^n/\cW_1 \rightarrow \bdC^n/\cW_1$),
and, with probability one, the minimal polynomial for any vector
in ~$\bdC^n/\cW_1$.
~To derive the induced map ~$\tilde{A}$,
~let ~$\{\,\bdq_{j\pls 1},\cdots,\bdq_n\,\}$~ be an orthonormal basis for
~$\spn{\bdx, A\bdx, \cdots, A^{j\mns 1}\bdx}^\perp$.
~Then
~$\Big\{ \, \bdx, A\bdx, \cdots, A^{j\mns 1}\bdx, \bdq_{j\pls 1},\cdots,\bdq_n
\Big\}$
~forms a basis for ~$\bdC^n$, ~and by writing ~$\tilde{\bdv} =
\bdv + \cW_1 \in \bdC^n/\cW_1$~ for any vector ~$\bdv \in \bdC^n$,
~$\{\,\tilde{\bdq}_{j\pls 1},\cdots,\tilde{\bdq}_n\,\}$~ forms a basis for
~$\bdC^n/\cW_1$.
~For the matrix representation of ~$\tilde{A}\,:\,\bdC^n/\cW_1 \rightarrow
\bdC^n/\cW_1$, ~let
\begin{eqnarray}
A\bdq_i & = &
c_{1i}\bdx + c_{2i}A\bdx + \cdots + c_{ji} A^{j\mns 1}\bdx +
c_{j\pls 1,i} \bdq_{j\pls 1}+\cdots+ c_{ni} \bdq_n
~~\mbox{\ \ for $i>j$} \nonumber \\
& = & \bdb_i + \big(c_{j\pls 1,i} \bdq_{j\pls 1}+\cdots+ c_{ni} \bdq_n\big)
\label{biq}
\end{eqnarray}
with ~$\bdb_i = c_{1i}\bdx + c_{2i}A\bdx + \cdots + c_{ji}A^{j\mns 1}\bdx
\in \cW_1$.
~It follows that
\[ \tilde{A}\tilde{\bdq}_i ~~=~~ \widetilde{A\bdq}_i
~~=~~ c_{j\pls 1,i} \tilde{\bdq}_{j\pls 1}+\cdots+c_{ni} \tilde{\bdq}_n 
\]
and the ~$(n-j)\times (n-j)$~ matrix
\[ \left[ \mbox{\scriptsize
$\begin{array}{lcl} c_{j\pls 1,j\pls 1} & \cdots & c_{j\pls 1,n} \\
\;\; \vdots & \ddots & \;\; \vdots \\
c_{n,j\pls 1} & \cdots & c_{n,n} \end{array}$} \right]
\]
becomes the matrix representation of the linear transformation
~$\tilde{A}\,:\,\bdC^n/\cW_1 \rightarrow \bdC^n/\cW_1$
~with respect to the basis
~$\{\,\tilde{\bdq}_{j\pls 1},\cdots,\tilde{\bdq}_n\,\}$.
~Meanwhile, by (\ref{biq}), ~$c_{li} = \bdq_l^\h A \bdq_i$, ~for
~$l,i = j\lpls 1,\cdots,n$.
~With the matrix representation of ~$\tilde{A}\,:\,\bdC^n/\cW_1
\rightarrow \bdC^n/\cW_1$~ available, we may find the minimal
polynomial ~$p_2(t)$ ~for ~$\bdC^n/\cW_1$~ (with respect to ~$\tilde{A}$)
by following the same procedure that produces minimal polynomial
~$p_1(t)$~ for ~$\cW_1$~ (with respect to ~$A$).
~For instance, using generically chosen ~$\bdy \in \cW_1$, ~write
~$\bdy = y_1\bdx + y_2 A\bdx + \cdots + y_j A^{j\mns 1}\bdx +
y_{j\pls 1} \bdq_{j\pls 1} + \cdots + y_n \bdq_n$
~and consider ~$\tilde{\bdy} = (y_{j\pls 1},\cdots,y_n)^\top \in
\bdC^n/\cW_1$~ ($\simeq \bdC^{n\mns j}$).
~Checking the sequence of Krylov subspaces
~$\spn{\tilde{\bdy},\tilde{A}\tilde{\bdy}}$, ~$\spn{\tilde{\bdy},\tilde{A}
\tilde{\bdy}, \tilde{A}^2\tilde{\bdy}}$, ~$\cdots$
consecutively. 
~Let
~$\spn{\tilde{\bdy},\tilde{A}\tilde{\bdy},\cdots\tilde{A}^l\tilde{\bdy}}$~
be the first one with its dimension less than the number of generating vectors.
~That is, the relation ~$ d_0 \tilde{\bdy} + d_1\tilde{A}\tilde{\bdy} + \cdots +
d_l \tilde{A}^l\tilde{\bdy} = \bdo$ ~with ~$d_l \ne 0$ ~exists, and
polynomial
~$\tilde{p}_2(t) = t^l + \frac{d_{l\mns 1}}{d_l} t^{l\mns 1} + \cdots +
\frac{d_{0}}{d_l}$
~becomes the minimal polynomial for ~$\tilde{\bdy}$.
~With probability one, it is the minimal polynomial for ~$\bdC^n/\cW_1$~
(with respect to ~$\tilde{A}$).
~Therefore ~$p_2(t) = \tilde{p}_2(t)$.

Notice that the linear independence of
~$\{\, \tilde{\bdy}, \tilde{A}\tilde{\bdy}, \cdots,
\tilde{A}^{l\mns 1}\tilde{\bdy}\,\}$~ in
~$\bdC^n/\cW_1$~ implies the linear independence of
~$\{\, \tilde{\bdy}, \tilde{A}\tilde{\bdy}, \cdots,
\tilde{A}^{l\mns 1}\tilde{\bdy}\,\}$~ in ~$\bdC^n$.
~Thus ~$\cW_2 = \spn{\bdy, A\bdy, \cdots, A^{l\mns 1}\bdy}$
~and
~$\cW_1 \oplus \cW_2 = \spn{\bdx, A\bdx, \cdots, A^{j\mns 1}\bdx,
\bdy, A\bdy, \cdots A^{l\mns 1}\bdy }$.
~In general,
~$\bdC^n/(\cW_1\oplus\cdots\oplus\cW_{m\mns 1}) \simeq
\cW_m\oplus\cdots\oplus\cW_k$, ~for ~$m = 2,\cdots,k$,
~so the same process may be continued to find the minimal polynomial
~$p_i(t)$~ for ~$\cW_i$, $i=3,\cdots,k$.

%In summary, ~let ~$\{\,\la_1,\cdots,\la_l\,\} = \eig{A}$. 
%%
%~For each (distinct) eigenvalue ~$\la_i \in \eig{A}$
%~with Segre characteristic ~$\{n_{i1}\geq  n_{i2} \ge \cdots \}$, 
%~the sequence of the minimal polynomials are
%\[ p_j(t) ~~=~~ (t-\la_1)^{n_{1j}} (t-\la_2)^{n_{2j}}\cdots (t-\la_l)^{n_{lj}},
%\;\;\;\; j = 1, 2, \cdots, \infty,
%\]
%where the nontrivial part of the sequence
%~$\big\{\,p_j(t)\,\big\}_{j=1}^\infty$~ is finite.

\vspace{-5mm}
\subsection{The minimal polynomial via Hessenberg reduction} \label{s:minpoly}
\vspace{-5mm}

In the process elaborated in the last section (\S \ref{sec:mp}), a crucial
step for
finding minimal polynomials is the determination of the dimensions of
the Krylov subspaces spanned by vector sets 
~$\{\bdx, A\bdx, \cdots, A^{j\mns 1} \bdx\}$ ~for ~$j=1,2,\cdots$.
~However, the condition of the Krylov matrix \label{'krylov'}
~$K(A,\bdx,j) \,\equiv\,
\blb \bdx, A\bdx, \cdots, A^{j\mns 1}\bdx \brb$ ~deteriorates
when ~$j$~ increases, making the rank decision difficult.
~A more reliable method is developed below to decide the dimension
of ~$\spn{\bdx, A\bdx, \cdots, A^{j\mns 1}\bdx}$~ accurately without
the explicit calculation of the Krylov matrices.

Computing eigenvalues of a matrix ~$A \in \bdC^{n\tms n}$~ starts with
the Hessenberg reduction \cite[p.344]{gvl}
\begin{equation} \label{hesred}
 Q^\h A Q ~~=~~ H ~~=~~ \blb \bdh_1,\cdots,\bdh_n \brb,
\;\;\; \mbox{with \ } Q^\h Q = I. \end{equation}
Let ~$\bdq_1, \cdots, \bdq_n $~ be the column vectors of ~$Q$~ in
(\ref{hesred}).
~Then
\begin{eqnarray*}
\lefteqn{
Q^\h \blb \bdq_1, A \bdq_1, \cdots, A^{j\mns 1} \bdq_1 \brb ~~=~~
   \blb \bde_1, Q^\h A \bdq_1, \cdots, Q^\h A^{j\mns 1} \bdq_1 \brb} \\
 & = & \blb \bde_1, (Q^\h A Q)Q^\h \bdq_1, \cdots, (Q^\h A^{j\mns 1} Q)Q^\h 
\bdq_1
\brb  ~~=~~  \blb \bde_1, H \bde_1, \cdots, H^{j\mns 1} \bde_1 \brb ~~=~~ R_{j}.
\end{eqnarray*}
Here ~$\bde_1 = [1,0,\cdots,0]^\top$.
~Clearly, ~$R_{j}$~ is an ~$n\times j$~ upper triangular matrix.
~Therefore
\begin{equation} \label{qrj}
 K(A,\bdq_1,j) ~~=~~
\blb \bdq_1, A \bdq_1, \cdots, A^{j\mns 1} \bdq_1 \brb ~~=~~ Q R_{j}
\end{equation}
is a QR decomposition of the Krylov matrix
~$K(A,\bdq_1,j) = \blb \bdq_1, A \bdq_1, \cdots, A^{j\mns 1} \bdq_1 \brb$.
~Furthermore, if ~$K(A,\bdq_1,j)$~ is of full rank, then
the first ~$j$ ~columns ~$\bdq_1,\cdots,\bdq_{j}$ ~of
~$Q$ ~form an orthonormal basis for
the Krylov subspace ~$\cR{K(A,\bdq_1,j)}$. 

Taking (\ref{qrj}) into account for computing the minimal polynomial
via Krylov matrices ~$K(A,\bdv, j)$ ~for ~$j=1,2,\cdots,$~
using randomly chosen unit vector~$\bdv$, 
~the Hessenberg reduction matrix ~$Q$~
in (\ref{hesred}) must have ~$\bdv$~ as its first column.
~This can be achieved by a modified Hessenberg reduction
\begin{equation} \label{heskry}
\left\{ \begin{array}{l}
\mbox{Find the Householder matrix ~$T$ ~such that 
~$T^\h \bdv = \bde_1$} \\
\mbox{for ~$j=1,2,\cdots$ ~do} \\
~\left\lfloor \begin{array}{l} 
\mbox{Hessenberg reduction \cite[\S 7.4.3]{gvl} step ~$j$ ~on ~$T^\h A T$: 
~Obtaining ~$P_j$ ~so} \\
\mbox{that the first ~$j$ column block ~$[\bdh_1,\cdots,\bdh_j]$ ~of 
~$(P_1 \cdots P_j)^\h (T^\h A T) (P_1 \cdots P_j)$} \\
\mbox{is upper-Hessenberg}
\end{array} \right.
\end{array} \right.
\end{equation}
Since ~$\bdv = T \bde_1$, ~the first column of ~$T$~
is the same as ~$\bdv$.
~The subsequent Hessenberg reduction steps of ~$T^\h A T$~ with unitary
transformations ~$P_1\cdots P_j$ ~does not change its first ~$j-1$ 
~columns. 
~Consequently the first ~$j$-column block ~$[\bdq_1,\cdots,\bdq_j]$~ of
~$T P_1\cdots P_k$ ~stay the same for ~$k\ge j$ ~with ~$\bdq_1 = \bdv$. 
~Upon completing (\ref{heskry}) for ~$j$ ~up to ~$n$, 
~we obtain the Hessenberg matrix ~$Q^\h A Q = H$~ with a
specified first column ~$\bdq_1 = \bdv$~ in ~$Q = TP_1 \dots P_{n}$. 

When the Krylov matrix ~$K(A,\bdv,j)$~ is of full rank, then
~$\cR{K(A,\bdv,j)}  =  \cR{[\bdq_1,\cdots,\bdq_j]}$ ~and
~$\cR{K(A,\bdv,j\lpls 1)} = \cR{\bdq_1,A\blb \bdq_1,\cdots,\bdq_j \brb}$.
~Thus, the rank of ~$K(A,\bdv,j\lpls 1)$ ~can be decided by
finding the numerical rank of
~$\big[ \bdq_1, A[ \bdq_1,\cdots,\bdq_j] \big]$ ~during the process 
(\ref{heskry}).
~Moreover, ~$AQ = QH$~ implies
~$A\blb \bdq_1,\cdots,\bdq_j\brb
= Q \blb \bdh_1,\cdots,\bdh_j \brb$
~where ~$\bdh_1,\cdots, \bdh_n$~ are columns of ~$H$.
~Consequently
\[ \big[ \bdq_1, A[ \bdq_1,\cdots,\bdq_j] \big]
~=~ \big[ \bdq_1, Q [ \bdh_1,\cdots,\bdh_j ] \big]
~=~ Q \big[ \bde_1, \bdh_1,\cdots,\bdh_j \big].
\]
Therefore, the numerical rank of
~$\big[ \bdq_1, A[ \bdq_1,\cdots,\bdq_j] \big]$~
is the same as the {\em upper-triangular} matrix
~$\big[ \bde_1, \bdh_1,\cdots,\bdh_j \big]$.
~We summarize this result in the following proposition.

\vspace{-4mm}
\begin{prop}
~For ~$A \in \bdC^{n \tms n}$,
~let ~$Q$~ be the unitary transformation matrix whose first column is
parallel to ~$\bdv \in \bdC^n$~ such that ~$Q^\h A Q = H = \blb \bdh_1,\cdots, 
\bdh_n \brb$~ is upper-Hessenberg.
~Assume ~$j > 0$~ is the smallest integer
for which Krylov matrix ~{\em $K(A,\bdv,j\lpls 1)$}~ is rank-deficient, then
~$\rank{K(A,\bdv,i)} = \rank{\blb \bde_1, \bdh_1,\cdots,\bdh_{i\mns 1} \brb}$
~for ~$i = 2, \cdots, j$.
\end{prop}
\vspace{-4mm}

When Krylov matrices ~$K(A,\bdv,i)$ ~for ~$i \le j$ ~are of full rank,
the matrix ~$\blb \bde_1, \bdh_1,\cdots,\bdh_j \brb$~ is rank-deficient 
in exact sense if and only if the diagonal entry ~$h_{j\mns 1,j}$ ~is zero
since the matrix ~$\blb \bde_1, \bdh_1,\cdots,\bdh_j \brb$~ is upper-triangular.
~In numerical computation,
however, an upper-triangular matrix can be numerically rank deficient
even though its diagonal entries are not noticeably small, e.g.,
the Kahan matrix \cite[p.260]{gvl}.
~That is, ~$h_{j\mns 1,j}$~ is usually small but not near zero 
for
~$\blb \bde_1, \bdh_1,\cdots,\bdh_j \brb$~ to be numerically rank
deficient.
~Therefore we must apply the inverse iteration (\ref{nvec})
to determine whether
~$\blb \bde_1,\bdh_1,\cdots,\bdh_j \brb$~ is rank
deficient in approximate sense.

When the first index ~$j$ ~is encountered with
~$\blb \bde_1,\bdh_1,\cdots,\bdh_j \brb$~ being numerically 
rank-deficient, we can further refine the Hessenberg reduction 
and minimize the magnitude of the entry ~$h_{j\pls 1, j}$ ~of ~$H$ ~since
\begin{equation} \label{hesodeq}
A \blb \bdq_1,\cdots,\bdq_j \brb - \blb \bdq_1,\cdots,\bdq_j \brb 
\blb \hat{\bdh}_1,\cdots,\hat{\bdh}_j\brb ~=~ 
h_{j\pls 1,j} \bdq_{j\pls 1} \end{equation}
should be zero, here ~$\hat{\bdh}_i \in \bdC^j$ ~is the first ~$j$-entry
subvector of ~$\bdh_i$ ~for ~$i=1,\cdots,j$.
~As a result, the least squares solution to the overdetermined system
\begin{equation} \label{hessys}
\left\{ \begin{array}{rcl}
A \blb \bdq_1,\cdots,\bdq_j \brb - \blb \bdq_1,\cdots,\bdq_j \brb 
\blb \hat{\bdh}_1,\cdots,\hat{\bdh}_j\brb &=& \bdo \\
\bdq_i^\h\, \blb \bdc_1,\cdots,\bdc_{i\mns 1}, \bdc_i\brb -
\blb 0, \cdots, 0, 1 \brb &=& \bdo, ~~~\mbox{for}~~ i=1,\cdots,j, \\
\bdq_1 - \bdv &=& \bdo \\
h_{il} & = & 0, ~~~\mbox{for}~~ i> l\lpls 1
\end{array} \right.
\end{equation}
minimizes ~$h_{j\pls 1, j}$. 
~Here ~$\bdv$ ~is the predetermined random vector and ~$\bdc_1,\cdots,\bdc_j$
~are constant vectors.
~Let ~$\bdf(\bdq_1,\cdots,\bdq_j,\bdh_1,\cdots,\bdh_j)$ ~be the 
vector mapping that represents the left side of the system (\ref{hessys})
~and ~$J(\bdq_1,\cdots,\bdq_j,\bdh_1,\cdots,\bdh_j)$ be its Jacobian. 
~The following proposition ensures the local convergence of the Gauss-Newton
iteration in solving 
~$\bdf(\bdq_1,\cdots,\bdq_j,\bdh_1,\cdots,\bdh_j) = \bdo$ ~for the least 
squares solution.

\vspace{-4mm}
\begin{prop}
~Let ~$A \in \bdC^{n\tms n}$ ~and 
~$\bdf(\bdq_1,\cdots,\bdq_j,\bdh_1,\cdots,\bdh_j) = \bdo$ ~be the vector
form of the system {\em (\ref{hessys})}. 
~Assume ~$\bdq_1,\cdots,\bdq_j,\bdh_1,\cdots,\bdh_j$ ~satisfies 
{\em (\ref{hessys})} ~with ~$h_{i\pls 1,i} \ne 0$ ~for
~{\em $i=1,\cdots,j\lmns 1$}.  
~Then the Jacobian of ~$\bdf(\bdq_1,\cdots,\bdq_j,\bdh_1,\cdots,\bdh_j)$ ~is 
injective.
\end{prop}
\vspace{-4mm}

\prf
~Differentiating the system (\ref{hessys}), let 
matrix ~$\blb \bdz_1,\cdots,\bdz_j \brb \in \bdC^{n\tms j}$ 
~and upper-Hessenberg matrix 
~$\blb \bdg_1,\cdots,\bdg_j \brb \in \bdC^{j\tms j}$
~satisfy
\begin{equation} \label{hess1}
\left\{ \begin{array}{rcl}
A\blb \bdz_1,\cdots,\bdz_j \brb  - \blb \bdz_1,\cdots,\bdz_j \brb 
\blb \hat\bdh_1,\cdots,\hat\bdh_j \brb & = & \blb \bdq_1,\cdots,\bdq_j \brb 
\blb \bdg_1,\cdots,\bdg_j \brb \\
\blb \bdc_1,\cdots,\bdc_{i\mns 1}, \bdc_i\brb^\h \bdz_i 
&=& \bdo, ~~~\mbox{for}~~ i=1,\cdots,j, \\
\bdz_1 ~=~ \bdo, ~~~~
g_{il} & = & 0, ~~~\mbox{for}~~ i> l\lpls 1
\end{array} \right.
\end{equation}
Using an induction, we have ~$\bdz_1=\bdo$ ~and assume 
~$\bdz_1 = \cdots = \bdz_k = \bdo$. 
~The equation ~$A\bdz_k = \sum_{i=1}^{k\pls 1}(h_{ik}\bdz_i +
g_{ik}\bdq_i)$ ~becomes
~$h_{k\pls 1,k} \bdz_{k\pls 1} + \sum_{i=1}^{k\pls 1} g_{ik}\bdq_i = \bdo$.
~For ~$i=1,\cdots,k$, ~we have ~$g_{ik} = 0$ ~from 
~$\bdc_i^\h \bdz_{k\pls 1} = 0$, ~$\bdc_i^\h \bdq_i=1$, 
~$\bdc_i^\h \bdq_l = 0$ ~for ~$l=i\lpls 1, \cdots, k\lpls 1$. 
~Also, ~$\bdc_{k\pls 1}^\h \bdz_{k\pls 1} = 0$ ~and 
~$\bdc_{k\pls 1}^\h \bdq_{k\pls 1} = 1$ ~lead to ~$g_{k\pls 1,k} = 0$ ~and
~$\bdz_{k\pls 1} = \bdo$ ~since ~$h_{k\pls 1,k} \ne 0$.
~Thus ~$\bdz_i = \bdo$ ~and ~$\bdg_i = \bdo$ ~for ~$i=1,\cdots,j$. 
~Namely, ~$J(\bdq_1,\cdots,\bdq_j,\bdh_1,\cdots,\bdh_j)$ ~is injective.
\qed

When ~$\blb \bde_1,\bdh_1,\cdots,\bdh_j\brb$ ~is numerically 
rank-deficient, we set ~$\blb \bdc_1,\cdots,\bdc_j\brb = 
\blb \bdq_1,\cdots,\bdq_j\brb$
~and the initial iterate 
~$\bdz^{(0)} = \blb \bdq_1^\h,\cdots,\bdq_j^\h,\bdh_1^\h,\cdots,\bdh_j^\h \
\brb^\h$
~for the Gauss-Newton iteration
\begin{eqnarray} \label{hesit}
\lefteqn{\bdz^{(i\pls 1)} ~=~ \bdz^{(i)} - 
J\big(\bdz^{(i)}\big)^+\bdf\big(\bdz^{(i)}\big)} \\
&& 
~~~\mbox{with}~~ \bdz^{(i)} = \blb \big(\bdq_1^{(i)}\big)^\h, \cdots, 
\big(\bdq_j^{(i)}\big)^\h, 
\big(\bdh_j^{(i)}\big)^\h, \cdots, \big(\bdh_j^{(i)}\big)^\h \brb^\h ,   
~~i=0,1,\cdots
\nonumber
\end{eqnarray}
that refines the (partial) Hessenberg reduction 
~$A\blb \bdq_1,\cdots,\bdq_j\brb = \blb \bdq_1,\cdots,\bdq_j\brb
\blb \hat\bdh_1,\cdots,\hat\bdh_j\brb$ 
~and minimize the magnitude of the residual
~$h_{j\pls 1,j} \bdq_{j\pls 1}$ ~that approaches zero during the iterative 
refinement. 

Overwrite ~$\bdq_1,\cdots,\bdq_j,\bdh_1,\cdots,\bdh_j$ ~with the
terminating iterate of (\ref{hesit}) and 
~$U_1\, R_1$ ~be the QR decomposition of 
~$\blb \bdq_1,\cdots,\bdq_j \brb$.
~Then 
\[ 
 U_1^\h A U_1 = \left[ \begin{array}{cc} H_1 & * \\ O & A_1 \end{array}
\right] 
\]
with ~$H_1$ ~being an upper-Hessenberg matrix whose characteristic 
polynomial is the first minimal polynomial ~$p_1$ ~of ~$A$.
~By the argument in \S~\ref{sec:mp}, the second minimal polynomial
~$p_2$ ~of ~$A$ ~is the (first) minimal polynomial of ~$A_1$. 
~Therefore we can continue the same Hessenberg reduction-refinement strategy 
on ~$A_1$ ~recursively and
obtain a reduced-Hessenberg form 
\begin{equation} \label{redhes}
U^\h A U = \left[ \begin{array}{ccc} H_1 & \cdots & * \\
& \ddots & \vdots \\  & & H_\ell \end{array} \right]
\end{equation}
where each ~$H_i$ ~is an irreducible upper-Hessenberg matrix whose 
characteristic polynomial is the ~$i$-th minimial polynomial ~$p_i$ ~of ~$A$
~for ~$i=1,\cdots,\ell$.

The first minimal polynomial
~$p_1(t) = p_0 + p_1 t + \cdots + p_j t^j$
~and its coefficient vector
~$\bdp \equiv (p_0,p_1,\cdots,p_j)^\top$~ satisfies
~$K(A,\bdv,j\lpls 1) \bdp = \bdo$.
~From (\ref{qrj}), ~$\bdv = \bdq_1$, ~and
~$K(A,\bdv,j)$~ being full rank, we have ~$K(A,\bdv,j) = Q R_{j} =
Q \mbox{\scriptsize $\left( \begin{array}{c} \hat{R}_{j} \\ O \end{array}
\right) $}$ ~where ~$\hat{R}_j$~ is a ~$j\times j$ ~upper triangular matrix and
\begin{eqnarray*}
K(A,\bdv,j\lpls 1) & = & \Big[ \bdv, ~A K(A,\bdv,j) \Big]
~=~ \Big[ \bdv, ~AQR_{j} \Big] ~=~
\Big[ \bdv, ~Q\blb \bdh_1,\cdots,\bdh_{j}\brb \hat{R}_{j} \Big] \\
& = & \Big[ Q\bde_1,Q\blb \bdh_1,\cdots,\bdh_{j}\brb \hat{R}_{j} \Big]
~=~
Q [\bde_1,H_1] \mbox{\scriptsize $\left[ \begin{array}{cc} 1 & \\ & \hat{R}_{j}
\end{array} \right]$}.
\end{eqnarray*}
In general, to find the coefficient vector ~$\bdp_i$ ~of the ~$i$-th 
minimal polynomial ~$p_i$ ~for ~$i=1,\cdots,\ell$, 
~we first solve
~$[\bde_1,H_i]\, \bdz = \bdo$ ~for ~$\bdz$ ~and write
\[
\bdp_i ~=~ \mbox{\scriptsize $
\left[ \begin{array}{c} \al \\ \bdu \end{array} \right]$},
\mbox{\ \ and \ \ }
\bdz ~=~ \mbox{\scriptsize
$\left[ \begin{array}{c} \al \\ \bdv \end{array} \right]$}.
\]
Then solve
\begin{equation} \label{mpcoef}
 \left[ \begin{array}{cc} 1 & \\ & \hat{R}_{j}
\end{array} \right]
\left[ \begin{array}{c} \al \\ \bdu \end{array} \right]
~~=~~ \left[ \begin{array}{c} \al \\ \bdv \end{array} \right].
\end{equation}

\vspace{-5mm}
\begin{itemize} \parskip-0.5mm
\item[] {\bf Algorithm} {\sc MinimalPolynomials}
\item[] {\tt Input}: ~$A \in \bdC^{n\tms n}$, ~numerical rank threshold 
~$\theta > 0$
\begin{itemize}
\item \hspace{-2mm} Initialize ~$\ell = 0$
\item \hspace{-2mm} While ~$n>0$ ~do
\item[] ~~$\left\lfloor \begin{array}{l}
\mbox{Set unit vector ~$\bdv$ ~at random. ~Apply (\ref{heskry}) until
~$\ranka{\theta}{[\bde_1,\bdh_1,\cdots,\bdh_j]} = j$} \\
\mbox{Update ~$\ell = \ell + 1$} \\
\mbox{Apply the Gauss-Newton iteration (\ref{hesit}) to refine
~$\bdq_1,\cdots,\bdq_j,\bdh_1,\cdots,\bdh_j$} \\
\mbox{Obtain the QR decomposition ~$U_\ell R_\ell = [\bdq_1,\cdots,\bdq_j]$}\\
\mbox{Obtain ~$U_\ell^\h A U_\ell = \left[ 
\mbox{\scriptsize $\begin{array}{cc} H_\ell & * \\ O & A_\ell
\end{array}$} \right]$.  ~Overwrite ~$A$ ~with ~$A_\ell$ ~and
update ~$n = n\lmns j$} \\
\mbox{Solve ~$[\bde_1,H_\ell]\,\bdz = \bdo$ ~for ~$\bdz$ ~and
construct ~$p_\ell$ ~by solving (\ref{mpcoef})} 
\end{array} \right.$
\end{itemize}
\item[] {\tt Output}: ~minimal polynomials ~$p_1, \cdots, p_\ell$
\end{itemize}
\vspace{-4mm}

The sequence of minimal polynomials ~$p_1(t)$, $\cdots$, ~$p_\ell(t)$~ 
produced by Algorithm {\sc MinimalPolynomials} are in the form
\[
 p_i(t)  ~~=~~  (t-\la_1)^{n_{1i}} \cdots (t-\la_l)^{n_{li}},
\;\;\; i = 1, 2 \cdots, \ell
\]
where ~$\{n_{j1} \geq n_{j2} \geq  \cdots\}$~ is the Segre characteristic
of ~$A$~ associated with ~$\la_j$ ~for $j=1,\cdots,l$.
~Although the process is recursive, there is practically no loss of
accuracy from ~$A_i$~ to ~$A_{i\pls 1}$~ since ~$A_{i\pls 1}$~ is extracted 
as a submatrix of ~$A_i$~ during which only unitary similarity transformations
are involved. 

For each ~$p_i(t)$, ~Algorithm {\sc MultRoot} in \cite{zeng_multroot,zeng-05}
is applied to calculate the multiplicity structure
~$[\,m_{i1},\cdots,m_{i\sg_i}\,]$~ and corresponding approximate roots
~$t_{i1}, \cdots,t_{i\sg_i}$, ~obtaining the Jordan structure of matrix ~$A$.

{\bf Remark.}
~The modified Hessenberg reduction (\ref{heskry}) is in fact the
Arnoldi process \cite[p. 172-179]{SaadBook} with Householder 
orthogonalization, which is the most reliable version of the Arnoldi
method.
~We improve its robustness even further with a novel iterative refinement
step (\ref{hesit}).
~There are less reliable versions of the Arnoldi iteration 
(see, e.g. \cite[p.303]{dem-book}\cite[p.499]{gvl}) based on Gram-Schmidt
orthogonalization that
may be applied to construct unitary bases 
for the Krylov subspaces.
~A method of finding minimal polynomials can alternatively be based 
on those versions of the Arnoldi algorithm.
~We choose the modified Hessenberg reduction and Gauss-Newton refinement
to ensure the highest possible accuracy. \qed

\vspace{-5mm}
\subsection{Minimal polynomials and matrix bundle stratification}
\vspace{-5mm}

In fact, the process of applying Algorithm {\sc MinimalDegree} on
matrix sequence ~$A_1$, ~$A_2$, ~$\cdots$ ~inherently calculates the
Segre characteristics associated with the matrix bundle of the highest
codimension.
~Suppose ~$A \in \bdC^{n\tms n}$ ~belongs to matrix bundle ~${\cal B}$ ~defined
by Segre characteristics
~$\{n_{j1} \ge n_{j2} \ge \cdots\}$, ~for ~$j=1,2,\cdots,k$.
~As explained in \S \ref{sec:ajcfnotion}, bundle ~${\cal B}$ ~is
imbedded in the closure of a lower codimension matrix bundle, say
~$\tilde{\cal B}$, ~in a hierarchy of bundle stratification.
~Our algorithm actually identifies the highest codimension
bundle ~${\cal B}$ ~because of the covering relationship established in
\cite{eek2}.

For minimal polynomials ~$p_1,p_2,\cdots$ ~of ~$A$, ~let ~$d_i \,=\,
\deg(p_i) \,=\, \sum_{j=1}^k n_{ji}$ ~for ~$i=1,\cdots,\infty$.
~The integer sequence ~$\{d_1 \ge d_2 \ge \cdots \}$ ~forms a partition
of ~$n$.
~Let ~$\{\tilde{d}_1 \ge \tilde{d_2} \ge \cdots \}$ ~be the similarly
constructed sequence of minimal polynomial degrees associated with
bundle ~$\tilde{\cal B}$ ~where
~$\overline{\tilde{\cal B}} \supseteq {\cal B}$.

\vspace{-3mm}
\begin{lemma} ~Suppose ~${\cal B}$ ~and ~$\tilde{\cal B}$ ~are two bundles
of ~$n\times n$ ~matrices with ~$\overline{\tilde{\cal B}} \supseteq {\cal B}$
~and a matrix on ~${\cal B}$ ~has at least as many distinct
eigenvalues as a matrix on ~$\tilde{\cal B}$.
~Let ~$d = \{d_1\ge d_2 \ge \cdots \}$ ~and
~$\tilde{d} = \{\tilde{d}_1 \ge \tilde{d}_2 \ge \cdots \}$
~be the degree sequences of minimal polynomials associated with
~${\cal B}$ ~and ~$\tilde{\cal B}$ ~respectively.
~Then ~$d$ ~and ~$\tilde{d}$ ~as partitions of ~$n$ ~satisfy the
dominant ordering relationship ~$ \tilde{d} \ge d$, ~namely
\begin{equation} \label{mpdom}
 \tilde{d}_1 + \tilde{d}_2 + \cdots + \tilde{d}_j ~~\ge ~~
d_1 + d_2 + \cdots + d_j ~~~~\mbox{for each}~~ j = 1, 2, \cdots.
\end{equation}
\end{lemma}
\vspace{-3mm}

\prf ~By \cite[Theorem 2.6]{eek2},
~$\overline{\tilde{\cal B}} \supseteq {\cal B}$
~if and only if it is possible to coalesce eigenvalues and apply the
dominance ordering coin moves to the Segre characteristics which defines 
bundle
~$\tilde{\cal B}$ ~to reach those of ~${\cal B}$.
~If ~${\cal B}$ ~is obtained by one
dominance coin move from one Segre characteristic
~$\tilde{\nu} = \{\tilde{n}_{j1} \ge \tilde{n}_{j2} \ge \cdots\}$ ~to
~$\nu = \{n_{j1} \ge n_{j2} \ge \cdots\}$ ~with other Segre characteristics
unchanged, ~then ~$\tilde{\nu} > \nu$ ~and
therefore (\ref{mpdom}) holds.

Similarly, assume ~${\cal B}$ ~is obtained by coalescing two eigenvalues on
~$\tilde{\cal B}$ ~with their Wyre characteristics combined as a union of
sets, or equivalently, their Segre characteristics 
~$\{\tilde{n}_{i1},\tilde{n}_{i2},\cdots\}$ ~and
~$\{\tilde{n}_{j1},\tilde{n}_{j2},\cdots\}$ ~combined in a componentwise
sum ~$\{\tilde{n}_{i1} + \tilde{n}_{j1}, ~\tilde{n}_{i2}+\tilde{n}_{j2}, 
\cdots \}$ 
and other Segre characteristics unchanged (see also \cite[Lemma 2.5]{eek2}).
~Actually the equalities in (\ref{mpdom}) hold
in this case since the degree ~$d_k$ ~is the sum 
~$n_{1k}+n_{2k} + \cdots$ ~of the ~$k$-th components in the Segre 
characteristics.

Since (\ref{mpdom}) is valid for every single dominant coin move and every
coalesce of eigenvalues,
it holds for a sequence of such manipulations
of Segre characteristics from ~$\tilde{\cal B}$ ~to ~${\cal B}$.
\hfill {\LARGE $\Box$} %\end{proof}

Because of (\ref{mpdom}), we have either ~$d = \tilde{d}$ ~or ~$d > \tilde{d}$.
~If ~$\tilde{d} > d$, there is an ~$l > 0$ ~such that
~$\tilde{d}_1 = d_1, ~\cdots, ~\tilde{d}_{l-1} = d_{l-1}$ ~and
~$\tilde{d}_l > d_l$.
~Algorithm {\sc MinimalDegree} applying on ~$A_l$ ~stops
at ~$d_l$ ~instead of ~$\tilde{d}_l$, ~since the search goes through
degree 1, 2, $\cdots$ ~and ~$d_l$ ~precedes ~$\tilde{d}_l$.
~Consequently, the highest codimension bundle ~${\cal B}$ ~is identified
before ~$\tilde{\cal B}$ ~with proper rank calculation.
~If ~$d = \tilde{d}$, ~the degrees of minimal polynomials associated with
~${\cal B}$ ~is the same as those of ~$\tilde{\cal B}$.
~For a similar reason,
Algorithm {\sc MultRoot} \cite{zeng_multroot} ~extracts
the highest codimension multiplicity structure that leads to ~${\cal B}$
~rather than ~$\tilde{B}$.
~Consequently, the highest codimension bundle ~${\cal B}$ ~is identified
before ~$\tilde{\cal B}$ ~with proper rank calculation.
~Our computing experiment is consistent with this observation.

\vspace{-5mm}
\section{The overall algorithm and numerical results} \label{s:overall}

\vspace{-5mm}
\subsection{The overall algorithm}
\vspace{-5mm}

Our overall algorithm for computing the numerical Jordan Canonical Form
of given matrix ~$A \in \bdC^{n\tms n}$ ~can now be summarized as follows.

\begin{itemize}
\item[] \hspace{-8mm} {\sc Stage I:} ~{\bf Computing the Jordan Structure}
\begin{itemize}
\item[] \hspace{-8mm} {\sc Step 1} ~~{\bf Francis QR.}
~Apply Francis QR algorithm to obtain a Schur decomposition
~$A \;=\; Q\,T\,Q^\h$
~and approximate eigenvalues ~$\la_1,\cdots,\la_n$.
\item[] \hspace{-8mm} {\sc Step 2} ~~{\bf Deflation.}
For each well-conditioned simple eigenvalue ~$\la_j$, ~apply the
deflation method
in \cite{bai-dem} to swap ~$\la_j$~ downward along the diagonal of ~$T$~
to reach
\[ A \;=\; U\,\mbox{\scriptsize $
\left[\begin{array}{cc} B & D \\ & C \end{array} \right]$} \, U^\h \]
where ~$\eig{C}$~ consists of all the well-conditioned eigenvalues of ~$A$.
\item[] \hspace{-8mm} {\sc Step 3} ~~{\bf Jordan structure.} ~Apply the method
in \S \ref{s:minpoly} to calculate the Segre characteristics
of ~$B$~ and initial estimates of the distinct eigenvalues.
\end{itemize}
\item[] \hspace{-8mm} {\bf \sc Stage II:}
~{\bf Computing the staircase/Jordan decompositions}
\newline 
There is an option here to select either the unitary staircase decompsition
~$A = USU^\h$~ or the Jordan decomposition ~$A = XJX^{\mns 1}$.
\begin{itemize}
\item[] \hspace{-8mm} {\sc Step 4(a)} 
~{\bf To compute the staircase decomposition.}
~~For each distinct eigenvalue, apply Algorithm
{\sc InitialEigentriplet} for an initial eigentriplet using
the Segre characteristic and initial eigenvalue approximation computed in the
previous step.
~Then iteratively refine the eigentriplet by Algorithm
{\sc EigentripletRefine}.
~Continue this process to reach a unitary staircase decomposition
~$A \;=\; U\,S\,U^\h$~ ultimately.
\item[] \hspace{-8mm} {\sc Step 4(b)} 
~{\bf To compute the Jordan decomposition.}
~For each distinct eigenvalue ~$\la_j$, ~$j=1,\cdots,k$~ with Segre
characteristic and the initial approximate determined in Step 3
above, apply the precess described in \S \ref{s:ctrip} to compute a
unitary-staircase eigentriplet ~$(\la_j,U_j,S_j)$.
~Then apply the Kublanovskaya algorithm to obtain the local Jordan decomposition
 ~$\la_j I + S_j \,=\, G_j S_j G_j^{\mns 1}$.
~Consequently, the Jordan decomposition ~$A \;=\; X\,J\,X^{\mns 1}$
~with ~$X = \blb U_1 G_1, \cdots, U_k G_k \brb$ ~is constructed.
\end{itemize}
\end{itemize}

As mentioned before, {\sc Stage II} can be considered a stand-alone
algorithm for computing the staircase/Jordan form from given Weyr/Segre
characteristics and initial eigenvalue approximation.
~It can be used in conjunction with other approaches
where the Jordan structure is determined by alternative means.

There are four control parameters that can be adjusted to improve the
results:

\vspace{-3mm}
\begin{enumerate} \parskip-1mm
\item {\em The deflation threshold ~$\dl$}: ~If a simple eigenvalue
has a condition number less than ~$\dl$, it will be deflated.
~The default value for ~$\dl$ ~is ~$1000$.
\item {\em The gap threshold in rank decision ~$\gamma$}: ~In determining
the rank deficiency of ~$H_j$~ in Algorithm {\sc MinimalDegree},
we calculate the smallest singular value of each ~$H_j$.
If the ratio of the smallest singular values of ~$H_j$~ and ~$H_{j\mns 1}$~
is less than ~$\gamma$, then ~$H_j$~ is considered rank deficient.
~The default value for ~$\gamma$ ~is ~$10^{-4}$.
\item {\em The residual tolerance ~$\tau$~ for {\sc MultRoot}}: ~The
residual tolerance required by {\sc MultRoot}. See \cite{zeng-05} for
details.
\item {\em The residual tolerance ~$\rho$~ for eigentriplet refinement}:
~It is used to stop the iteration in refining the eigentriplet.
The default value ~$\rho = 10^{\mns 8}$.
\end{enumerate}

\vspace{-5mm}
\subsection{Numerical results}
\vspace{-5mm}

We made a Matlab implementation {\sc NumJCF} of our algorithm for computing the
Jordan decomposition.
~It has been tested in comparison with the Matlab version of JNF
\cite{kagstrom-ruhe-jnf} on a large number of matrices, including
classical examples in the literature.
~Our experiment is carried out on a Dell Optiplex GX 270 personal computer with
Intel Pentium 4 CPU, 2.66 GHz and 1.5 GB RAM.
~For a computed Jordan decomposition ~$A X = XJ$~ of matrix ~$A$, ~the
residual ~$\rho ~=~ \|\,A\,X-X\,J\,\|_F\big/\|A\|_F $
~is used as one of the measures for the accuracy%
%\footnote{The Matlab code of our algorithm and test scripts for all the 
%numerical examples can be 
%accessed at {\tt http://www.neiu.edu/$\sim$zzeng/numjcf.htm}
%}.

\vspace{-4mm}
\begin{example} \em ~~Let
\[ A_4 \,=\,\left[ \mbox{\tiny $\begin{array}{rrrrrr}
2\,r-5-s & -r+3\,s-2\,t & 20-2\,s+2\,t & 15-2\,s+2\,t & 10 & -5+s-t\\
2\,r-5-2\,s & -r-15+6\,s-4\,t & 50-4\,s+4\,t & 40-4\,s+4\,t & 20 & -15+2\,s-2\,t\\
0 & -10-2\,s+2\,t & 10+4\,s-3\,t & 10+3\,s-3\,t & s-t & -5-s+t\\
2\,r-5-2\,s & -r-10+8\,s-7\,t & 50-8\,s+8\,t & 40-7\,s+8\,t & 25-s+t & -15+3\,s-3\,t\\
-2\,r+5+2\,s & r+25-6\,s+5\,t & -65+4\,s-4\,t & -55+4\,s-4\,t & -25+t & 25-2\,s+2\,t\\
0 & -5 & 10 & 10 & 5 & -5+t
\end{array}$} \right] \in \bdR^{6\tms 6}.
\]
We compare our method with the conventional symbolic computation
on the exact matrix.
~The exact eigenvalues are ~$r$, ~$s$~ and ~$t$~ with Segre characteristics
~$\{1\}$, ~$\{2\}$~ and ~$\{3\}$, ~respectively.
~For ~$r=\sqrt{2}$, ~$s=\sqrt{3}$~ and ~$t=\sqrt{5}$, ~it takes Maple 10 nearly
two hours (7172 seconds) to find the Jordan Canonical Form, while both JNF and
{\sc NumJCF} complete the computation instantly.
~On a similarly constructed matrix of size ~$10\times 10$, ~Maple does not
finish the computation in 8 hours and Mathematica runs out of memory.
\end{example}
\vspace{-4mm}

\begin{table}[ht] \small \begin{center}
\begin{tabular}{|l|l|l|l|c|} \hline
 & \multicolumn{3}{c|}{computed eigenvalues} & residual \\
 & \multicolumn{3}{c|}{(with correct digits in boldface and Jordan block sizes
in braces)} & $\rho$ \\ \hline
JNF     & \footnotesize {\bf 1.41421356}3 \{1\} &
\footnotesize {\bf 1.73205080}9 \{2\}&
\footnotesize {\bf 2.23606797}5 \{3\} & 3.38e-013 \\
{\sc NumJCF}  & \footnotesize {\bf 1.4142135623731}1 \{1\} &
\footnotesize {\bf 1.73205080757}4 \{2\}&
\footnotesize {\bf 2.2360679774997}1 \{3\}& 1.01e-016 \\ \hline
\end{tabular}
\caption{\footnotesize
Eigenvalues and residuals computed by {\sc JNF} and {\sc NumJCF}}
\label{tbl:1}
\vspace{-3mm}
\end{center} \end{table}

Approximating ~$\sqrt{2}, \;\sqrt{3}$~ and ~$\sqrt{5}$~ in machine precision
~$\approx 2.2\ltms 10^{-16}$,
~both JNF and {\sc NumJCF} correctly identify the Jordan
structure, whereas our {\sc NumJCF} obtained the eigenvalues with
$3 \sim 6$ more correct digits than JNF along with smaller residual as shown in
Table \ref{tbl:1}.

\begin{example} \em
~~This is a classic test matrix that is widely used in
eigenvalue computing experiment
\cite{chatelin-fraysse,kagstrom-ruhe,lzw,ruhe-70-bit,varah67}:
\[ A_5 ~~=~~ \left[ \mbox{\tiny $\begin{array}{rrrrrrrrrr}
  1  &  1  &  1  & -2  &  1  & -1  &  2  & -2  &  4  & -3 \\
 -1  &  2  &  3  & -4  &  2  & -2  &  4  & -4  &  8  & -6 \\
 -1  &  0  &  5  & -5  &  3  & -3  &  6  & -6  & 12  & -9 \\
 -1  &  0  &  3  & -4  &  4  & -4  &  8  & -8  & 16 & -12 \\
 -1  &  0  &  3  & -6  &  5  & -4  & 10 & -10  & 20 & -15 \\
 -1  &  0  &  3  & -6  &  2  & -2  & 12 & -12  & 24 & -18 \\
 -1  &  0  &  3  & -6  &  2  & -5  & 15 & -13  & 28 & -21 \\
 -1  &  0  &  3  & -6  &  2  & -5  & 12 & -11  & 32 & -24 \\
 -1  &  0  &  3  & -6  &  2  & -5  & 12 & -14  & 37 & -26 \\
 -1  &  0  &  3  & -6  &  2  & -5  & 12 & -14  & 36 & -25
\end{array}$} \right] \;\;\;\;\;
\mbox{ \begin{tabular}{|c||cc|} \hline
$\eig{A_5}$ & \multicolumn{2}{c|}{Segre ch.} \\ \hline
 1 & 1 & \\ 2 & 3 & 2 \\ 3 & 2 &  2 \\ \hline \end{tabular} }
\]
Using the default parameters, both {\sc JNF} and our {\sc NumJCF} easily
obtained
the accurate Jordan decomposition.
\end{example}

\scriptsize
\begin{center} \begin{tabular}{|l|c|c|c|} \hline
\multicolumn{4}{|c|}{\normalsize Computing results for eigentriplets of
~$A_5$~ } \\ \hline
&   eigenvalue & Segre characteristic & residual \\ \hline
&1.0000000000000002 &  \{1, 0\} &          \\
{\sc JNF}
&2.0000000000000001 &  \{3, 2\} & 1.41e-15 \\
&3.0000000000000002 &  \{2, 2\} &          \\ \hline
&0.9999999999999995 &  \{1, 0\} &          \\
{\sc NumJCF}
&2.0000000000000000 &  \{3, 2\} & 1.40e-16 \\
&3.0000000000000003 &  \{2, 2\} &          \\ \hline
\end{tabular} \end{center}\normalsize
\vspace{-3mm}

Both {\sc JNF} and {\sc NumJCF} obtain
similarly accurate results
on classical matrices such as those in \cite[pp.
192-196]{chatelin-fraysse}.
~We choose to omit them and concentrate on the cases in which our
{\sc NumJCF} significantly improves the robustness and accuracy in
comparison with {\sc JNF}.

\begin{example} \em
~~This is a series of test matrices with a parameter ~$t$.
\begin{equation}
A(t) ~~=~~
\mbox{\scriptsize $\left[ \begin{array}{rrrrrrrrrr}
     t&    2\pls t&  \mns t&     \mns 2& \mns 1\mns 3t&     \mns 2&      2&   \mns 1\pls t&    \mns t&    0 \\
   1\mns t& \mns 1\mns 3t& 2t&    2\pls t&  2\pls 6t&    2\pls t&   \mns 3\mns t&    1\mns t& 1\pls 2t&    1  \\
   2t&   \mns 4t&   2&    4t&      t&    3t&   \mns 2t&      t&     0&    0  \\
  \mns 1\pls t&   \mns 7t& 2t&  1\pls 4t& 1\pls 10t& \mns 1\pls 4t& \mns 1\mns 5t&      0& 1\pls 3t&  1\pls t  \\
   3t&   \mns 4t&   0&    4t&    3\pls t&    4t&  1\mns 3t&    2t&     0&   \mns 1  \\
 2\mns 3t& \mns 4\pls 5t&   0&  4\mns 4t&  1\mns 5t&  6\mns 4t& \mns 2\pls 6t&    1\mns t&    \mns t& \mns 2t  \\
\mns 3\pls 4t&  2\mns 3t&  \mns t& \mns 2\pls 4t& \mns 1\mns 3t& \mns 2\pls 4t&  6\mns 2t& \mns 1\pls 4t&    \mns t& \mns 1\mns t  \\
   4t&   \mns 5t&   0&    5t&      t&    5t&  1\mns 4t&  3\pls 3t&     0&   \mns 1  \\
\mns 2\mns 3t& \mns 2\pls 2t&   t&   \mns 3t&  1\pls 2t&   \mns 3t& \mns 3\pls 2t&   \mns 2t&   4\pls t&    3  \\
\mns 3\pls 4t&  2\mns 3t&  \mns t& \mns 2\pls 4t& \mns 1\mns 3t& \mns 2\pls 4t&  3\mns 2t& \mns 1\pls 4t&    \mns t&  2\mns t
\end{array} \right]$} \in \bdR^{10\tms 10}  \label{exam11}
\end{equation}
For every ~$t > 0$, ~matrix ~$A(t)$~ has the same Jordan Canonical Form ~$J$~
consisting of two eigenvalues ~$\la_1 = 2$~ and ~$\la_2 = 3$~ with Segre
characteristics ~$\{3,1\}$~ and ~$\{4,2\}$~ respectively.
~Let the Jordan decomposition be ~$A(t) = X(t)\,J\,X(t)^{\mns 1}$.
~When ~$t$~ increases, the condition number ~$\|X(t)\|_2\|X(t)^{\mns 1}\|_2$~
of ~$X(t)$~ increases rapidly.
~This example tests the accuracy and robustness of numerical Jordan Canonical
Form finders under the increasing condition number of ~$X(t)$.
~As shown in Table \ref{tb:exam11}, our algorithm maintains high backward
accuracy, forward accuracy, and structure correctness, while the results of
{\sc JNF} deteriorate as ~$t$~ increases.
~Starting from ~$t=5$, {\sc JNF} outputs incorrect Jordan structure.
~When ~$t\ge 23$, {\sc JNF} outputs only one ~$10\times 10$~ Jordan block.
~Our {\sc NumJCF} continues to produce accurate results.
\end{example}

\begin{table}[ht]
\begin{center} \scriptsize
\begin{tabular}{|c|l|l|l|c|l|} \hline
\multicolumn{2}{|c|}{} &  & &  & \\
\multicolumn{2}{|c|}{} & eigenvalues & Segre ch. & backward error & $\|X(t)\|_2\|X(t)^{\mns 1}\|_2$ \\
\multicolumn{2}{|c|}{} &  & &  & \\  \hline \hline
 &    & 2.00000000000001 & \ \ \ \ {3,1} & & \\
      & {\sc JNF}  & 2.99999999999999 & \ \ \ \ {4,2} & 6.10e-015 &
\\ \cline{2-5}
 $t=1$     & & 2.00000000000000 & \ \ \ \ {3,1} & & \ \ 1113.9 \\
      & {\sc NumJCF} & 3.00000000000000 & \ \ \ \ {4,2} & 1.11e-015 & \\ \hline \hline
 &    & 2.000000000002 & \ \ \ \ {3,1} & & \\
      &  {\sc JNF} & 2.999999999998 & \ \ \ \ {4,2} & 7.40e-014 & \\ \cline{2-5}
  $t=2$    & & 2.00000000000000 & \ \ \ \ {3,1} & & \ \ 28894.5 \\
      &  {\sc NumJCF} & 3.00000000000000 & \ \ \ \ {4,2} & 4.87e-016 & \\ \hline \hline
 &    & 1.99999999987 & \ \ \ \ {3,1} &  &\\
      &  {\sc JNF} & 3.00000000009 & \ \ \ \ {4,2} & 8.30e-012 & \\ \cline{2-5}
  $t=4$    & & 2.00000000000000 & \ \ \ \ {3,1} & & \ \ 1658396.6 \\
      &  {\sc NumJCF}  & 2.99999999999999 & \ \ \ \ {4,2} & 5.65e-016 & \\ \hline \hline
 &    & 2.0000000006 & \ \ \ \ {4} & & \\
      &  {\sc JNF} & 2.9999999996 & \ \ \ \ {5,1} & 1.36e-011 & \\ \cline{2-5}
 $t=5$     & & 2.00000000000001 & \ \ \ \ {3,1} & & \ \ 5655648.5 \\
      &  {\sc NumJCF} & 2.99999999999999 & \ \ \ \ {4,2} & 7.60e-016 & \\ \hline \hline
 &  & 2.0000001  & \ \ \ \ {4} & & \\
       &  {\sc JNF}  & 2.99999992 & \ \ \ \ {5,1} & 9.29e-010 & \\ \cline{2-5}
 $t=10$     & & 2.00000000000003 & \ \ \ \ {3,1} & & \ \ 297244917.4 \\
      &  {\sc NumJCF}       & 2.99999999999998 & \ \ \ \ {4,2} & 6.94e-016 & \\ \hline \hline
 & {\sc JNF}    & 2.6  & \ \ \ \ {10} & 1.27e-004& \\  \cline{2-5}
 $t=25$     & & 1.99999999999992 & \ \ \ \ {3,1} & & \ \ 60948418207.9 \\
      &  {\sc NumJCF} & 2.99999999999998 & \ \ \ \ {4,2} & 8.58e-016 & \\ \hline \hline
\end{tabular}
\end{center}
\vspace{-3mm}
\caption{\footnotesize
Comparison between {\sc JNF} and {\sc NumJCF} on matrix ~$A(t)$~ in
(\ref{exam11})}
\label{tb:exam11}
\vspace{-3mm}
\end{table}

\vspace{-3mm}
\begin{example} \label{e:1000} \em
~~The matrices in the literature on computing Jordan Canonical Forms
are usually not larger than ~$10\times 10$.
~We construct a ~$100 \times 100$~ real matrix
\[  A \;=\; X \left[ 
\mbox{\scriptsize $\begin{array}{cc} J & \\ & B \end{array}$}
\right] X^{\mns 1},
\]
where ~$J$~ is the Jordan Canonical Form of eigenvalues ~$\la_1 = 1$~ and
~$\la_2 = 2$~ with Segre characteristics ~$\{5,4,3,1\}$~ and ~$\{4,2,2\}$~
respectively, ~$B \in \bdR^{80 \tms 80}$~ and ~$X \in \bdR^{100 \tms 100}$~
are random matrices with entries uniformly distributed in ~$[-1,1]$.
~There are 80 simple eigenvalues randomly scattered around the two
multiple eigenvalues.
~This example is designed to show that our code {\sc NumJCF} may be more
reliable
than the approach of grouping the eigenvalue clusters in the process of
identifying a multiple eigenvalue and determining the Jordan structure.
~We generate 1000 such matrices ~$A$~ with fixed ~$J$~ and randomly chosen ~$B$~
as well as ~$X$.
~For each matrix ~$A$, ~we run {\sc JNF} and {\sc NumJCF}
twice and the results are shown in Table \ref{1000mat}.
\end{example}

\begin{table}[ht]
\begin{center}
\begin{tabular}{|l|c|c|c|} \hline
& \% of failures & \% of failures & \% of failures \\
& on both run & on first run & on second run \\ \hline
{\sc JNF} & 41.9\% & 41.9\% & 41.9\% \\
{\sc NumJCF} & 0.1\% & 4.5\% & 4.6\% \\ \hline
\end{tabular}
\end{center}
\vspace{-3mm}
\caption{\footnotesize
Results for Example \ref{e:1000} on 1000 matrices.} \label{1000mat}
\end{table}

Notice that there are several steps in our algorithm which require parameters
generated at random.
~Consequently, failures are rarely repeated (0.1\% in this case) in the
subsequent runs of {\sc NumJCF}.
The code {\sc JNF} appears to be deterministic and always repeats the same 
results.
~On the other hand, failures are verifiable in our algorithm from the residuals
staircase condition numbers.
~One may simply run the code second time when the first run fails.

\end{document}